\documentclass[MSNbibl,number,citesort,seceqn,dvips]{arxbj}
\usepackage{graphicx}

\aid{0}
\volume{20}
\issue{3}
\pubyear{2014}
\firstpage{1532}
\lastpage{1559}
\doi{10.3150/13-BEJ532} 

\makeatletter
\newcommand{\rrvert}{\vert}
\newcommand{\llvert}{\vert}
\def\cal{\mathcal}
\newtheorem{thmm}{Theorem}[section]
\newremark{rem}{Remark}[section]
\newproclaim{assumption}{Assumption}[section]
\newproclaim{condition}{Condition}[section]
\newtheorem{proposition}{Proposition}[section]
\newtheorem{corol}{Corollary}[section]
\newtheorem{lemma}{Lemma}[section]
\newremark{example}{Example}[section]

\newcommand{\field}[1]{\mathbb{#1}}
\newcommand{\R}{\field{R}}

\newcommand{\p}{\field{P}}
\newcommand{\N}{\field{N}}
\newcommand{\Z}{\field{Z}}

\newcommand{\E}{\field{E}}

\def\argmin{\mathop{\operatorname{argmin}}}

\makeatother

\begin{document}
\begin{frontmatter}

\title{Asymptotics of nonparametric L-1 regression models with
dependent data}
\runtitle{Asymptotics of nonparametric L-1 regression}

\begin{aug}
\author[a]{\inits{Z.}\fnms{Zhibiao} \snm{Zhao}\corref{}\thanksref{a,e1}\ead[label=e1,mark]{zuz13@stat.psu.edu}},
\author[b]{\inits{Y.}\fnms{Ying} \snm{Wei}\thanksref{b,e2}\ead[label=e2,mark]{ying.wei@columbia.edu}}
\and
\author[a]{\inits{D.K.J.}\fnms{Dennis K.J.} \snm{Lin}\thanksref{a,e3}\ead[label=e3,mark]{DKL5@psu.edu}}
\runauthor{Z. Zhao, Y. Wei and D.K.J. Lin} 
\address[a]{Department of Statistics, Penn State University, University
Park, PA 16802, USA.\\
\printead{e1,e3}}

\address[b]{Department of Biostatistics, Columbia University, 722 West
168th St., New York, NY 10032, USA.
\printead{e2}}
\end{aug}

\received{\smonth{12} \syear{2011}}
\revised{\smonth{5} \syear{2013}}

%
\begin{abstract}
We investigate asymptotic properties of least-absolute-deviation
or median quantile estimates of the location and scale functions
in nonparametric regression models with dependent data from multiple
subjects. Under a general dependence
structure that allows for longitudinal data and some spatially
correlated data,
we establish uniform Bahadur representations for the
proposed median quantile estimates. The obtained Bahadur
representations provide deep insights
into the asymptotic behavior of the estimates. Our main theoretical
development is based on studying the modulus of continuity of
kernel weighted empirical process through a coupling argument.
Progesterone data is used for an illustration.
\end{abstract}

%

\end{frontmatter}

\section{Introduction}\label{sec:intro}
There is a vast literature on the nonparametric location-scale model
$Y=\mu(X)+s(X)e$,
where $X,Y$, and $e$ are the covariates, response, and error, respectively.
Given observations $\{(X_j,Y_j)\}_{j=1,\ldots,m}$, the latter model has
been studied under various settings of data structure.
In terms of the dependence structure, there are
independent data and time series data scenarios; in terms of the design
point $X$,
there are random-design and fixed-design $X_j=j/m$ settings.
In these settings, we usually assume that either $(X_j,Y_j)$ are
independent observations from
subjects $j=1,\ldots,m$, or $\{(X_j,Y_j)\}_{j=1,\ldots,m}$ is a
sequence of time series observations from the same subject.
We refer the reader
to Fan and Yao \cite{fan2003} and Li and Racine \cite{li2007} for an
extensive exposition of related works.

In this article, we are interested in
the following nonparametric location-scale model with serially
correlated data from multiple subjects:
%
%
\begin{equation}
\label{eq:model} Y_{i,j}=\mu(x_{i,j})+s(x_{i,j})e_{i,j},\qquad
1\le j\le m_i, 1\le i\le n,
\end{equation}
where, for each subject $i$, $\{(x_{i,j},Y_{i,j})\}_{j=1,\ldots,m_i}$ is
the sequence of covariates and responses, and $\{e_{i,j}\}_{j=1,\ldots
,m_i}$ is the corresponding
error process.
We study (\ref{eq:model}) under a general dependence
framework for $\{e_{i,j}\}_{j\in\N}$ that allows for both longitudinal
data and some spatially correlated data.
In typical longitudinal studies, $x_{i,j}$ represents measurement time or
covariates at time $j$, then it is reasonable to assume that
$\{e_{i,j}\}_{j\in\Z}$ is a causal time series, that is, the current
observation depends only on past but not future observations.
In other applications, however, measurements may be dependent on
both the
left and right neighboring measurements, especially when $x_{i,j}$
represents measurement location.
A good example of this type of data is the vertical density profile
data in
Walker and Wright \cite{wal2002}; see also Section~\ref{sec:exm} for
more details.
To accommodate this, we propose
a general error dependence structure, which can be viewed as an
extension of the one-sided causal structure in Wu \cite{wu2005} and
Dedecker and Prieur \cite{ded2005} to a two-sided noncausal setting.
The proposed
dependence framework allows for many linear and nonlinear
processes.

We are interested in
nonparametric estimation of the location function $\mu(\cdot)$ and the scale
function $s(\cdot)$. Least-squares based nonparametric methods have been
extensively studied for both time series data (Fan and Yao \cite
{fan2003}) and longitudinal data
(Hoover \textit{et al.} \cite{hoo1998},
Fan and Zhang \cite{fan2000},
Wu and Zhang \cite{wu2002},
Yao, M\"uller and Wang \cite{yao2005}). While they
perform well for Gaussian errors, least-squares based methods are
sensitive to extreme outliers, especially when
the errors have a heavy-tailed distribution. By contrast, robust
estimation methods impose
heavier penalty on far-deviated data points to reduce the
impact from extreme outliers. For example, median quantile
regression uses the absolute loss and the resultant estimator
is based on sample local median. Since Koenker and Bassett \cite{koe1978},
quantile regression has become
popular in parametric and nonparametric inferences and we refer
the reader to Yu, Lu and Stander \cite{yu2003} and Koenker \cite{koe2005}
for excellent expositions. Recently, He, Fu and Fung \cite{he2003},
Koenker \cite{koe2004} and Wang and Fygenson~\cite{wan2009}
applied quantile regression techniques to
parameter estimation of parametric longitudinal models,
He, Zhu and Fung \cite{he2002} studied median regression for
semiparametric longitudinal
models, and Wang, Zhu and Zhou \cite{wang2009}
studied inferences for a partially linear varying-coefficient
longitudinal model.
Here we focus on quantile regression based estimation for
the nonparametric model (\ref{eq:model}).

We aim to study the asymptotic properties, including uniform
Bahadur representations and asymptotic normalities, of the
least-absolute-deviation or median quantile estimates
for model~(\ref{eq:model}) under a general
dependence structure.
Nonparametric quantile regression estimation has been
studied mainly under either the i.i.d. setting (Bhattacharya and
Gangopadhyay~\cite{bha1990},
Chaudhuri \cite{cha1991},
Yu and Jones \cite{yu1998})
or the strong mixing setting (Truong and Stone~\cite{tru1992},
Honda~\cite{hon2000},
Cai \cite{cai2002}).
There are relatively scarce results on Bahadur representations of
conditional quantile estimates. Bhattacharya and Gangopadhyay \cite
{bha1990} and Chaudhuri \cite{cha1991}
obtained point-wise Bahadur representations for conditional quantile
estimation of i.i.d. data. For
mixing stationary processes, Honda \cite{hon2000} obtained point-wise and
uniform Bahadur representations of conditional quantile
estimates. For stationary random fields, Hallin, Lu and Yu \cite
{hal2009} obtained a point-wise Bahadur representation
for spatial quantile regression function under spatial mixing conditions.
Due to the nonstationarity and dependence
structure, it is clearly challenging to establish
Bahadur representations in the context of (\ref{eq:model}).

Our contribution here is mainly on the theoretical side. We establish
uniform Bahadur representations for the least-absolute-deviation
estimates of $\mu(\cdot)$ and
$\sigma(\cdot)$ in (\ref{eq:model}). To derive
the uniform Bahadur representations, the key ingredient is to
study the modulus of continuity of certain kernel weighted
empirical processes of the nonstationary observations $Y_{i,j}$
in (\ref{eq:model}). Empirical processes have been extensively
studied under various settings, including the i.i.d. setting (Shorack
and Wellner \cite{sho1986}),
linear processes (Ho and Hsing \cite{ho1996}),
strong mixing setting (Andrews and Pollard \cite{and1994},
Shao and Yu \cite{sha1996}),
and general causal stationary processes (Wu \cite{wu2008}). Using
a coupling argument to approximate the dependent process by
an $m$-dependent process with a diverging $m$, we study the modulus
of continuity of weighted empirical processes, and the
latter result serves as a key tool in establishing our uniform
Bahadur representations. These Bahadur representations provide
deep insights into the asymptotic behavior of the estimates,
and in particular they provide theoretical justification
for the profile control chart methodologies in Wei, Zhao and Lin \cite{wei2012}.
These technical treatments are also of
interest in other nonparametric problems involving dependent
data.

The article is organized as follows. In Section~\ref{sec:dependence}, we introduce the error dependence
structure with examples. In Section~\ref{sec:coupling}, we study
weighted empirical process
through a coupling argument.
Section~\ref{sec:main} contains
uniform Bahadur
representations and asymptotic normality.
Section~\ref{sec:app}
contains an illustration using progesterone data.
Possible extensions to spatial setting are discussed in Section~\ref{sec:dis}.
Proofs are provided in
Section~\ref{sec:proof}.

\section{Error dependence structure}\label{sec:dependence}
First, we introduce some notation used throughout
this article. For $a,b\in\R$, let $\lfloor a\rfloor$ be the
integer part of $a$, $a\vee
b=\max(a,b)$, and $a\wedge b=\min(a,b)$. For a random variable
$Z \in{\cal L}^q, q>0$, if $\|Z\|_q = [\E
(|Z|^q)]^{1/q} < \infty$. Let
${\cal C}^r ({\cal S})$ be the set of functions with bounded
derivatives up to order $r$ on a set ${\cal S}\subset\R$.

Assume that, for each $i$, the error process $\{e_{i,j}\}_{j\in
\N}$ in (\ref{eq:model}) is an independent copy from a
stationary process $\{e_j\}_{j\in\N}$ which has the
representation
%
%
\begin{equation}
\label{eq:error} e_j=G(\varepsilon_j,
\varepsilon_{j\pm1}, \varepsilon_{j\pm2},\ldots),
\end{equation}
where $\varepsilon_j,j\in\Z$, are i.i.d. random innovations, and
$G$ is a measurable function such that $e_j$ is well defined.
We can view (\ref{eq:error}) as an input-output system with
$(\varepsilon_j,\varepsilon_{j\pm1}, \varepsilon_{j\pm
2},\ldots), G$, and $e_j$ being, respectively, the input, filter,
and output. Wu \cite{wu2005} considered the causal time series case
that $e_j$ depends only on the past innovations
$\varepsilon_j,\varepsilon_{j-1},\ldots.$
In contrast, (\ref{eq:error}) allows for noncausal models and is particularly
useful for applications that do not have a time structure. For
example, if $x_{i,j}$ are locations, then the corresponding measurement
$y_{i,j}$ depends on
both the left and right neighboring measurements.

\begin{condition}\label{con:dc}
Let $\{\varepsilon'_j\}_{j\in\Z}$ be i.i.d. copies of
$\{\varepsilon_j\}_{j\in\Z}$. There exist constants $q>0$ and
$\rho\in(0,1)$ such that
%
%
\begin{eqnarray}
\label{eq:gmc} \bigl\|e_0-e_0(k)\bigr\|_q=\mathrm{O}\bigl(
\rho^k\bigr),\qquad \mbox{where } e_0(k) = G\bigl(
\varepsilon_0,\varepsilon_{\pm1},\ldots,
\varepsilon_{\pm
k}, \varepsilon'_{\pm(k+1)},
\varepsilon'_{\pm(k+2)},\ldots\bigr).\qquad
\end{eqnarray}
\end{condition}

In (\ref{eq:gmc}), $e_0(k)$ can be viewed as a coupling process
of $e_0$ with $\{\varepsilon_r\}_{|r|\ge k+1}$ replaced
by the i.i.d. copy $\{\varepsilon'_r\}_{|r|\ge k+1}$ while keeping
the nearest $2k+1$ innovations $\{\varepsilon_r\}_{|r|\le k}$.
In particular, if $e_0$ does not depend on
$\{\varepsilon_r\}_{|r|\ge k+1}$, then $e_0(k)=e_0$. Thus,
$\|e_0-e_0(k)\|_q$ quantifies the contribution of
$\{\varepsilon_r\}_{|r|\ge k+1}$ to $e_0$, and
(\ref{eq:gmc}) states that the contribution decays
exponentially in $k$. Shao and Wu~\cite{sha2007} and Dedecker and
Prieur \cite{ded2005} [cf. equation (4.2) therein]
considered one-sided causal version of (\ref{eq:gmc}) where $e_0$
depends only on
$\{\varepsilon_r\}_{r\le0}$.

Propositions \ref{pro:h}--\ref{pro:e} below indicate that, if $\{e_i\}$
satisfies (\ref{eq:gmc}), then its
properly transformed process also satisfies (\ref{eq:gmc}).

\begin{proposition}\label{pro:h}
For $0 < \varsigma\le1$ and $\upsilon\ge0$, define the
collection of functions $h$
%
%
\begin{equation}
\label{eq:Hfunc} {\cal H}(\varsigma, \upsilon) = \bigl\{h\dvt \bigl|h(x)-h
\bigl(x'\bigr)\bigr|\le c \bigl|x-x'\bigr|^\varsigma \bigl(1 +
|x|+\bigl|x'\bigr| \bigr)^\upsilon, x, x'\in\R\bigr\},
\end{equation}
where $c$ is a constant.
Suppose $\{e_j\}$ satisfies (\ref{eq:gmc}). Then the transformed
process $e^*_j=h(e_j)$ satisfies (\ref{eq:gmc}) with $(q,\rho)$
replaced by $q^*=q/(\varsigma+\upsilon)$ and $\rho^*=\rho^\varsigma$.
\end{proposition}

In (\ref{eq:Hfunc}), ${\cal H}(\varsigma,0)$ is the class of uniformly
H\"older-continuous functions with index $\varsigma$. If
$h(x)=|x|^b, b> 1$, then $h\in{\cal H}(1,b-1)$. Clearly, all functions
in ${\cal H}(\varsigma,0)$
are continuous. Interestingly, for noncontinuous transformations, the
conclusion may
still hold; see Proposition \ref{pro:e} below, where $\mathbf{1}$ is
the indicator function.

\begin{proposition}\label{pro:e}
Let $e_0$ have a bounded density. Suppose $\{e_j\}$ satisfies (\ref
{eq:gmc}). Then, for any given~$x$, $\{\mathbf{1}_{e_j\le x}\}$
satisfies (\ref{eq:gmc})
with $\rho$ replaced by $\rho^*=\rho^{1/(1+q)}$.
\end{proposition}

Propositions \ref{pro:h}--\ref{pro:e} along with the examples below
show that
the error structure (\ref{eq:error}) and Condition
\ref{con:dc} are sufficiently general to accommodate many
popular linear and nonlinear time series models and their properly
transformed processes.

\begin{example}[($\bolds{m}$-dependent sequence)]\label{exmp:mdep}
Assume that $e_j=G(\varepsilon_{j},\varepsilon_{j\pm
1},\ldots,\varepsilon_{j\pm m})$ for a measurable function $G$.
Then $e_j$ depends only on the nearest $2m+1$ innovations
$\varepsilon_j,\varepsilon_{j\pm1},\ldots,\varepsilon_{j\pm
m}$. Clearly, $\{e_j\}_{j\in\Z}$ form a $(2m+1)$-dependent
sequence, $\|e_0-e_0(k)\|_q=0$ for $k\ge m$, and
(\ref{eq:gmc}) trivially holds. If $m=0$, then $e_j$ are i.i.d.
random variables.
\end{example}

\begin{example}[(Noncausal linear processes)]\label{exmp:vlp}
Consider the noncausal linear process
$e_j=\sum^{\infty}_{r=-\infty} a_r \varepsilon_{j-r}$. If $\varepsilon
_j\in{\cal L}^q$ and $a_j=\mathrm{O}(\rho^{|j|})$, then
it is easy to see that
(\ref{eq:gmc}) holds.
\end{example}

\begin{example}[(Iterated random functions)]\label{exmp:irf}
Consider random variables $e_j$ defined
by
%
%
\begin{equation}
\label{eq:nar} e_j=R(e_{j-1},\ldots,e_{j-d};
\varepsilon_j),
\end{equation}
where $\varepsilon_j,j\in\Z$, are i.i.d. random innovations, and
$R$ is a random map. Many widely time series models are of form (\ref{eq:nar}),
including threshold autoregressive model
$e_j=a\max(e_{j-1}, 0) + b\min(e_{j-1}, 0) +
\varepsilon_j$, autoregressive conditional heteroscedastic
model $e_j=\varepsilon_j (a^2+b^2 e^2_{j-1})^{1/2}$,
random coefficient model $e_j=(a+b\varepsilon_j)
e_{j-1}+\varepsilon_j$, and exponential autoregressive model
$e_j=[ a+b\exp(-c e^2_{j-1})] e_{j-1} +\varepsilon_j$,
among others.
Suppose there exists $z_0$ such
that $R(z_0;\varepsilon_0)\in{\cal L}^q$ and there exist
constants $a_1,\ldots,a_d$ such that
\begin{eqnarray*}
\sum^d_{j=1} a_j <1\quad
\mbox{and}\quad\bigl \|R(z;\varepsilon_0)-R\bigl(z';
\varepsilon_0\bigr)\bigr\|^{1\wedge q}_q \le\sum
^d_{j=1} a_j \bigl|z_j-z'_j\bigr|^{1\wedge q}
\end{eqnarray*}
holds for all $z=(z_1,\ldots,z_d),z'=(z'_1,\ldots,z'_d)$. By Shao and
Wu \cite{sha2007},
(\ref{eq:gmc}) holds.
\end{example}

\subsection{Some examples}\label{sec:exm}

The imposed dependence structure and hence the developed results in
Sections~\ref{sec:coupling}--\ref{sec:main} below
are potentially applicable to a wide range of practical data types. We
briefly mention some below.

(\emph{Time series data}). In the special case of $n=1$, $m_1=m\to
\infty$ and $(x_{1,j},Y_{1,j},e_{1,j})=(x_j,Y_j,e_j)$
for a stationary time series $\{e_j\}$, (\ref{eq:model}) becomes
the usual nonparametric location-scale model $Y_j=\mu(x_j)+s(x_j) e_j$
with time series data. The latter model has been
extensively studied under both the random-design case and the
fixed-design case $x_j=j/n$. See Fan and Yao \cite{fan2003} for
an excellent introduction to various local least-squares based methods
under mixing settings.
Quantile regression based estimations have been studied in Truong and
Stone \cite{tru1992}, Honda \cite{hon2000}, and Cai \cite{cai2002} for
mixing processes.
Despite the popularity of mixing conditions, it is generally difficult
to verify mixing conditions even for linear processes.
For example, for the autoregressive model $X_i=\rho X_{i-1}+\varepsilon
_i, \rho\in(0,1/2]$, where $\varepsilon_i$ are i.i.d. Bernoulli random
variables
$\p(\varepsilon_i=1)=1-\p(\varepsilon_i=0)=q\in(0,1)$, the stationary
solution is not strong mixing (Andrews~\cite{and1984}).
By contrast, as shown above, the imposed Condition \ref{con:dc}
is easily verifiable for many linear and nonlinear time series models
and their proper transformations.

(\emph{Longitudinal data}). For each subject $i$, if $x_{i,j}$ is the
$j$th measurement time or the covariates at time $j$, $Y_{i,j}$ is the
corresponding response,
and $\{e_{i,j}\}_{j\in\N}$ is a stationary causal process [e.g.,
$e_j=G(\varepsilon_j,\varepsilon_{j-1}, \varepsilon_{j-2},\ldots)$ in
(\ref{eq:error}) depends only on the past], then (\ref{eq:model}) becomes
a typical longitudinal data setting. For example, Section~\ref{sec:app2} re-examines the well-studied
progesterone data using the proposed methods. Another popular
longitudinal data example is
the CD4 cell percentage in HIV infection from the Multicenter
AIDS Cohort Study. Based on least-squares methods, this data has been
studied previously in
Hoover \textit{et al.} \cite{hoo1998} and Fan and Zhang \cite{fan2000}.
We can examine how the response function (CD4 cell percentage)
varies with measurement time (age) using the proposed robust estimation
method in Section~\ref{sec:main}.

(\emph{Spatially correlated data}). In the vertical
density data of Walker and Wright \cite{wal2002},
manufacturers are concerned about engineered wood boards'
density, which determines fiberboard's overall quality.
For each board, densities are measured at various locations along a
designated vertical line.
In this example, measurements depend on both the left and right
neighboring measurements, and it is reasonable to
impose the dependence structure (\ref{eq:error}). See Wei, Zhao and Lin
\cite{wei2012} for a detailed analysis.
Also, as will be discussed in Section~\ref{sec:dis}, the two-sided
framework (\ref{eq:error}) can be extended to spatial lattice settings.
We point out that the structure in (\ref{eq:model}) and (\ref
{eq:error}) differs from the usual spatial model setting
in the sense that (\ref{eq:model}) allows for observations from
multiple independent subjects whereas the latter usually assumes that
all observations are spatially correlated (see, e.g., Hallin, Lu and Yu
\cite{hal2009} for
quantile regression of spatial data).

\section{Weighted empirical process}\label{sec:coupling}
In this section,
we study weighted empirical processes through a coupling argument.
Dependence is the main difficulty in extending results developed for
independent data to dependent data.
For mixing processes, the widely used large-block-small-block technique
partitions the data into asymptotically independent blocks.
Here, we adopt a coupling argument which copes well with the
dependence structure in Section~\ref{sec:dependence}.

We now illustrate the basic idea. By
(\ref{eq:error}), the error $e_{i,j}$ in (\ref{eq:model}) has
the representation
\[
e_{i,j}=G(\varepsilon_{i,j},\varepsilon_{i,j\pm1},
\varepsilon_{i,j\pm
2},\ldots)
\]
for i.i.d. innovations $\varepsilon_{i,j}, i,j\in\Z$. Thus, $\{e_{i,j}\}
_{j\in\Z}$ is a dependent series for each fixed $i$, whereas
$\{e_{i_1,j}\}_{j\in\Z}$ and $\{e_{i_2,j}\}_{j\in\Z}$ are two
independent series for $i_1\ne i_2$.
Let $\varepsilon'_{i,j,k}, i,j,k\in\Z$, be i.i.d.
copies of $\varepsilon_{i,j}$. For $k_n\in\N$, define the
coupling process of $e_{i,j}$ as
%
%
\begin{equation}
\label{eq:couple} e_{i,j}(k_n)= G\bigl(\varepsilon_{i,j},
\varepsilon_{i,j\pm1}, \ldots, \varepsilon_{i,j\pm k_n},
\varepsilon'_{i,j,j\pm(k_n+1)},\varepsilon '_{i,j,j\pm(k_n+2)},
\ldots\bigr)
\end{equation}
by replacing all but the nearest $2k_n+1$ innovations with i.i.d.
copies. We call $k_n$ the coupling lag. Clearly, $e_{i,j}(k_n)$
has the same distribution as $e_{i,j}$.

By construction, for each fixed $i$, $\{e_{i,j}(k_n)\}_{j\in\Z}$
form $(2k_n+1)$-dependent sequence in the sense that
$e_{i,j}(k_n)$ and $e_{i,j'}(k_n)$ are independent if $|j-j'|\ge
2k_n+1$. Consequently, for each fixed $i$ and $s$,
$\{e_{i,(j-1)(2k_n+1)+s}(k_n)\}_{j\in\Z}$ are i.i.d.
The latter property helps us reduce the dependent data to
an independent case. On the other hand, under (\ref{eq:gmc}),
$\|e_{i,j}-e_{i,j}(k_n)\|_q=\mathrm{O}(\rho^{k_n})$ is sufficiently small with
properly chosen $k_n$
and hence the coupling process is close
enough to the original one. Similarly, for $Y_{i,j}$ in (\ref
{eq:model}), define its coupling process:
%
%
\begin{equation}
\label{eq:yij2} \tilde{Y}_{i,j}=\mu(x_{i,j})+s(x_{i,j})e_{i,j}(k_n).
\end{equation}

First, we present a general result regarding the sum of functions of
the coupling process $\tilde{Y}_{i,j}$.
Let $\mathcal{V}_n$ be any finite set.
For real-valued functions $g_{i,j}(y,v), i,j\in\N$, defined on $\R
\times\mathcal{V}_n$
such that $\E[g_{i,j}(\tilde{Y}_{i,j},v)]=0$ for all $v\in\mathcal
{V}_n$, define
\[
H_n(v)=\sum^{n}_{i=1} \sum
^{m_i}_{j=1} g_{i,j}(
\tilde{Y}_{i,j},v),\qquad v\in \mathcal{V}_n.
\]
Throughout, let $N_n=m_1+\cdots+m_n$ be the total number of observations.

\begin{thmm}\label{thmm:cp}
Assume that the cardinality $|\mathcal{V}_n|$ of $\mathcal{V}_n$ and
the coupling lag $k_n$ grow no faster than a polynomial of $N_n$.
Further assume $|g_{i,j}(y,v)|\le c$ for a constant $c<\infty$, and for
some sequence $\chi_n$,
%
%
\begin{equation}
\label{eq:thmcpcon} \max_{v\in\mathcal{V}_n} \sum
^{n}_{i=1} \sum^{m_i}_{j=1}
\E \bigl[g^2_{i,j}(\tilde{Y}_{i,j},v)\bigr] \le
\chi_n.
\end{equation}
\begin{longlist}[(ii)]
\item[(i)] If $\chi_n=\mathrm{O}(1)$, then $\max_{v\in\mathcal{V}_n} |H_n(v)|=\mathrm{\textup{O}}_\mathrm
{p}(k_n \log N_n)$.

\item[(ii)] If $\sup_n \log N_n/\chi_n <\infty$, then $\max_{v\in\mathcal
{V}_n} |H_n(v)|=\mathrm{O}_\mathrm{p}[ k_n (\chi_n \log N_n)^{1/2}]$.
\end{longlist}
\end{thmm}

By Theorem \ref{thmm:cp}, the magnitude of $\max_{v\in\mathcal{V}_n}
|H_n(v)|$ increases with the coupling lag $k_n$.
Intuitively, as $k_n$ increases, there is stronger dependence in the
coupling process $\tilde{Y}_{i,j}$ and consequently
a larger bound for $H_n(v)$. Therefore, a small $k_n$ is preferred in
order to have a tight bound for $H_n(v)$. On the other hand,
a reasonably large $k_n$ is needed in order for the coupling process to
be sufficiently close to the original process.
Under (\ref{eq:gmc}), for $k_n=\mathrm{O}(\log N_n)$, the coupling process
converges to the original one
at a polynomial rate, and meanwhile the maximum bound in Theorem \ref
{thmm:cp} is optimal up to a logarithm factor. For example,\vspace*{1pt}
if $\chi_n=\mathrm{O}(1)$, then $\max_{v\in\mathcal{V}_n} |H_n(v)|=\mathrm{O}_\mathrm{p}[
(\log N_n)^2]$; if $\sup_n \log N_n/\chi_n <\infty$, \
then $\max_{v\in\mathcal{V}_n} |H_n(v)|=\mathrm{O}_\mathrm{p}\{ [\chi_n (\log
N_n)^3]^{1/2}\}$.

In what follows, we consider the special case of weighted empirical
process, which plays an essential role in quantile regression. Let
$\varpi_{i,j}(x)\ge0$ be nonrandom weights that may depend on~$x$. Consider the weighted empirical process
%
%
\begin{equation}
\label{eq:Fn} F_n(x,y) = \sum^{n}_{i=1}
\sum^{m_i}_{j=1} \varpi_{i,j}(x)
\mathbf{1}_{Y_{i,j}\le y}.
\end{equation}
To study $F_n(x,y)$, recall $\tilde{Y}_{i,j}$ in (\ref{eq:yij2}) and
define the
coupling empirical process
%
%
\begin{equation}
\label{eq:Fn2} \tilde{F}_n(x,y) = \sum^{n}_{i=1}
\sum^{m_i}_{j=1} \varpi_{i,j}(x)
\mathbf{1}_{\tilde{Y}_{i,j}\le y}.
\end{equation}
Under mild regularity conditions, Theorem \ref{lem:appr} below states
that $F_n(x,y)$ can be
uniformly approximated by $\tilde{F}_n(x,y)$ with proper choice of the
coupling lag $k_n$.

\begin{condition}\label{assump:r}
\textup{(i)} $\varpi_{i,j}(x)\le c$ uniformly for some constant $c<\infty$. \textup{(ii)}
$\mu(x_{i,j})$ is uniformly bounded. \textup{(iii)} $s(x_{i,j})>0$ is uniformly
bounded away from zero and infinity.
\end{condition}

\begin{thmm}\label{lem:appr}
Assume that Conditions \ref{con:dc} and \ref{assump:r} hold. In
(\ref{eq:couple}), let the coupling lag $k_n
=\lfloor\lambda\log N_n \rfloor$ for some $\lambda>(q+1)/[q\log
(1/\rho)]$,
where $N_n=m_1+\cdots+m_n$. Then
\[
\sup_{x,y\in\R} \bigl|F_n(x,y)-\tilde{F}_n(x,y)\bigr|
= \mathrm{O}_\mathrm{p}\bigl[(\log N_n)^2\bigr].
\]
\end{thmm}

To study asymptotic Bahadur representations of quantile regression estimates,
a key step is to study the modulus
of continuity of $F_n(x,y)$, defined by
%
%
\begin{eqnarray}
\label{eq:empdiff} D_n(\delta,x,y) = \bigl\{ F_n(x,y+
\delta)-\E\bigl[ F_n(x,y+\delta) \bigr] \bigr\} - \bigl\{
F_n(x,y) - \E\bigl[ F_n(x,y)\bigr] \bigr\}.
\end{eqnarray}
Intuitively, $D_n(\delta,x,y)$ measures the oscillation of the centered
empirical process $F_n(x,y)-\E[F_n(x,y)]$
in response to a small perturbation $\delta$ in $y$.

The dependence structure in Section~\ref{sec:dependence} along with the
coupling argument provides a convenient
framework to study $D_n(\delta,x,y)$. Recall $\tilde{F}_n(x,y)$ in (\ref
{eq:Fn2}).
For $D_n(\delta,x,y)$ in (\ref{eq:empdiff}), define its coupling process
%
%
\begin{eqnarray}
\label{eq:empdiff2} \tilde{D}_n(\delta,x,y) = \bigl\{
\tilde{F}_n(x,y+\delta)-\E\bigl[ \tilde {F}_n(x,y+\delta)
\bigr] \bigr\} - \bigl\{ \tilde{F}_n(x,y) - \E\bigl[
\tilde{F}_n(x,y)\bigr] \bigr\}.
\end{eqnarray}
Notice that $e_{i,j}(k_n)$ and $e_{i,j}$ have the same distribution, so
$\E[ F_n(x,y)]=\E[\tilde{F}_n(x,y)]$. By Theorem \ref{lem:appr}, it is
easy to see that,
uniformly over $x,y,\delta$,
%
%
\begin{eqnarray}
\label{eq:dd2}\bigl |D_n(\delta,x,y)-\tilde{D}_n(\delta,x,y)\bigr|
\le2 \sup_{x,y\in\R} \bigl| F_n(x,y) - \tilde{F}_n(x,y)\bigr|
= \mathrm{O}_\mathrm{p}\bigl[(\log N_n)^2\bigr].
\end{eqnarray}
Therefore, the asymptotic properties of $D_n(\delta,x,y)$ are similar
to that of $\tilde{D}_n(\delta,x,y)$, which can be
studied through Theorem \ref{thmm:cp}.

\begin{condition}\label{assump:r2}
\textup{(i)} $\varpi_{i,j}(\cdot)=0$ outside a common bounded interval for all
$i,j$. \textup{(ii)} There exist $\tau_n$ and $\phi_n$ such that
%
%
\begin{eqnarray}
\label{eq:conr2} \sup_{x\ne x'} \frac{|\varpi_{i,j}(x)-\varpi_{i,j}(x')|}{|x-x'|} \le
\tau_n \quad\mbox{and}\quad \sup_x \sum
^n_{i=1} \sum^{m_i}_{j=1}
\varpi^2_{i,j}(x) \le\phi_n.
\end{eqnarray}
\end{condition}

\begin{thmm} \label{pro:osc}
Assume that Conditions \ref{con:dc} and \ref{assump:r}--\ref{assump:r2}
hold. Further assume $\delta_n\to0$,\break
$\sup_{n}\log N_n/(\delta_n\phi_n)<\infty$, and that $1/\delta_n+\tau
_n$ grows no faster than a polynomial of
$N_n$.
Then
%
%
\begin{eqnarray}
\label{eq:osc} \sup_{|\delta|\le\delta_n, x,y\in\R} \bigl|D_n(\delta,x,y)\bigr|
=\mathrm{O}_\mathrm{p}\bigl\{ \bigl[\delta_n\phi_n (\log
N_n)^3\bigr]^{1/2}\bigr\}.
\end{eqnarray}
\end{thmm}

\section{Quantile regression and Bahadur representations}\label{sec:main}

For a random variable $Z$, denote by ${\cal Q}(Z)=\inf\{z\in\R,
\p(Z\le z)\ge1/2\}$ the median of $Z$, and similarly denote by
${\cal Q}(\cdot|\cdot)$ the conditional median operator. To
ensure identifiability of $\mu$ and $s$ in (\ref{eq:model}), without
loss of generality
we assume ${\cal Q}(e_{i,j})=0$
and ${\cal Q}(|e_{i,j}|)=1$.

Note that ${\cal Q}(Y_{i,j}|x_{i,j}=x)=\mu(x)$. Applying a
kernel localization technique, we propose the following
least-absolute-deviation or median quantile estimate of
$\mu(x)$:
%
%
\begin{equation}
\label{eq:muest} \hat\mu(x) = \argmin_{\theta}\sum
_{i=1}^{n} \sum_{j=1}^{m_i}
|Y_{i,j}-\theta| K_{b_n}(x_{i,j}-x),\qquad \mbox{where }
K_{b_n}(u)=K(u/b_n)
\end{equation}
for a nonnegative kernel function $K$ satisfying $\int_{\R}
K(u)=1$, and $b_n>0$ is a bandwidth. The estimate
$\hat\mu_{b_n}(x)$ pools together information across all
subjects, an appealing property especially when some subjects
have sparse observations.
By the Bahadur
representation in Theorem \ref{thmm:m} below, the bias term of
$\hat\mu(x)-\mu(x)$ is of order $\mathrm{O}(b^2_n)$. Following Wu and Zhao \cite
{wu2007}, we adopt a jackknife bias-correction technique. In
(\ref{eq:muest}), denote by $\hat\mu(x|b_n)$ and
$\hat\mu(x|\sqrt{2}b_n)$ the estimates of $\mu(x)$ using
bandwidth $b_n$ and $\sqrt{2}b_n$, respectively. The
bias-corrected jackknife estimator is
%
%
\begin{equation}
\label{eq:mbias} \tilde\mu_(x)=2\hat\mu(x|b_n)-\hat\mu(x|
\sqrt{2}b_n),
\end{equation}
which can remove the second-order bias term $\mathrm{O}(b^2_n)$ in $\hat\mu(x)$.

After estimating $\mu(\cdot)$, we estimate $s(\cdot)$ based on
residuals. Notice that ${\cal Q}(|e_{i,j}|)=1$ implies ${\cal
Q}(|Y_{i,j}-\mu(x)||x_{i,j}=x) =s(x)$. Therefore, we propose
%
%
\begin{equation}
\label{eq:sest} \hat{s}(x) = \argmin_{\theta}\sum
_{i=1}^{n} \sum_{j=1}^{m_i}
\bigl|\bigl |Y_{i,j}-\tilde\mu(x)\bigr| - \theta \bigr| K_{h_n}(x_{i,j}-x),
\end{equation}
where $h_n>0$ is another bandwidth, and $\tilde\mu(x)$ is the
bias-corrected jackknife estimator in (\ref{eq:mbias}). As in
(\ref{eq:mbias}), we adopt the following bias-corrected
jackknife estimator
%
%
\begin{equation}
\label{eq:vbias} \tilde{s}(x)=2\hat{s}(x|h_n)-\hat{s}(x|\sqrt{2}
h_n).
\end{equation}

\begin{rem}\label{rem:other}
By ${\cal
Q}(|Y_{i,j}-\mu(x_{i,j})||x_{i,j}=x) =s(x)$, an alternative estimator
of $s(x)$ is
%
%
\begin{equation}
\label{eq:sestot} \bar{s}(x) = \argmin_{\theta}\sum
_{i=1}^{n} \sum_{j=1}^{m_i}
\bigl| \bigl|Y_{i,j}-\tilde\mu(x_{i,j})\bigr| - \theta \bigr| K_{h_n}(x_{i,j}-x).
\end{equation}
The difference between (\ref{eq:sest}) and (\ref{eq:sestot}) is that
(\ref{eq:sest}) uses $\tilde\mu(x)$
whereas (\ref{eq:sestot}) uses $\tilde\mu(x_{i,j})$.
Since~$K$ has bounded support, only those $x_{i,j}$'s with
$|x_{i,j}-x|=\mathrm{O}(h_n)$ contribute to the summation in~(\ref{eq:sestot}). Thus, as $h_n\to0$ so that $x_{i,j}\to x$ and
$\tilde\mu(x_{i,j})\approx\tilde\mu(x)$, the two
estimators in (\ref{eq:sest}) and~(\ref{eq:sestot}) are expected to be
asymptotically close.
Our use of (\ref{eq:sest}) has some technical and computational
advantages. First,
the estimation error $\tilde\mu(x_{i,j})-\mu(x_{i,j})$ varies with
$(i,j)$, and thus it is technically more challenging
to study (\ref{eq:sestot}). Second, to implement (\ref{eq:sestot}), we
need to compute $\tilde\mu(\cdot)$ at each point $x_{i,j}$,
which requires solving a large number of optimization problems in (\ref
{eq:muest}) for a large data set.
By contrast, (\ref{eq:sest}) only requires estimation of $\tilde\mu
(\cdot)$ at those grid points $x$ at which
we wish to estimate $s(\cdot)$.
\end{rem}

To study asymptotic properties, we need
to introduce some regularity conditions. Throughout we write
${\cal S}_\epsilon([a,b])=[a+\epsilon,b-\epsilon]$ for an
arbitrarily fixed small $\epsilon>0$. Denote by $F_e$ and
$f_e=F'_e$ the distribution and density functions of $e_0$ in
(\ref{eq:error}), respectively. The assumption ${\cal
Q}(e_0)=0$ and ${\cal Q}(|e_0|)=1$ implies $F_e(0)=1/2$ and
$F_e(1)-F_e(-1)=1/2$.

\begin{condition}\label{con:reg}
Suppose that all measurement locations $x_{i,j}$ are
within an interval $[a,b]$, and order them as
$a=\tilde{x}_0 < \tilde{x}_1 <\cdots<\tilde{x}_{N_n}
<\tilde{x}_{N_n+1}= b$. Assume that
%
%
\begin{equation}
\label{eq:dense} \max_{0\le k\le N_n} \biggl|\tilde{x}_{k+1}-
\tilde{x}_k - \frac
{b-a}{N_n} \biggr|=\mathrm{O}\bigl(N^{-2}_n
\bigr),\qquad \mbox{where } N_n=m_1+\cdots+m_n.
\end{equation}
\end{condition}

Condition \ref{con:reg} requires that the pooled covariates $x_{i,j}$
should be approximately uniformly dense on
$[a,b]$, which is a natural condition since otherwise it would
be impossible to draw inferences for regions with very scarce observations.
Pooling all subjects together is an appealing
procedure to ensure this uniform denseness even though each
single subject may only contain sparse measurements.

In nonparametric regression problems, there are two typical settings on
the design points:
fixed-design and random-design points. For fixed-design case, it is
often assumed that the design points are equally spaced on some
interval. For example, for the
vertical density profile data of Walker and Wright \cite{wal2002},
the density was measured at equispaced points along a designated
vertical line of wood boards. Condition \ref{con:reg} can be viewed as
a generalization of the fixed-design points to allow for approximately
fixed-design points. For random-design case,
the design points are sampled from a distribution. For example,
assumption (a) in Appendix A of Fan and Zhang~\cite{fan2000} imposed
the random-design condition. In practice, both settings have different
range of applicability. For example, for daily or monthly temperature
series, the fixed-design setting may be appropriate; for
children's growth curve studies, it may be more reasonable to use the
random-design setting since the measurements are usually taken at
irregular time points.

\begin{rem}[(Asymptotic results under the random-design case)]\label{rem:random}
All our subsequent theoretical results are derived under the
approximate fixed-design setting in Condition \ref{con:reg}, but the
same argument also applies to the random-design case. Specifically,
assume that the design-points $\{x_{i,j}\}$ are random samples from a
density $f_X(\cdot)$ with support $[a,b]$ and that $x$ is an interior
point. Then, for the design-adaptive local linear median quantile
regression estimates, the subsequent
Theorems \ref{thmm:m}--\ref{thmm:v} and Corollaries \ref{con:reg}--\ref{con:lc} still hold with $(b-a)$
therein replaced by $1/f_X(x)$. In fact, given the i.i.d. structure of
$\{x_{i,j}\}$, the technical argument becomes much easier. For example,
to establish Lemma \ref{lem:1} (again, with $(b-a)$ therein replaced by
$1/f_X(x)$), elementary calculations can easily find the mean and
variance for the right-hand side of (\ref{eq:lem1a}). All other proofs can be
similarly modified and we omit the details.
\end{rem}

Conditions
\ref{con:lc}--\ref{con:regularity} below are standard assumptions in
nonparametric estimation.

\begin{condition}\label{con:lc}
$K$ is symmetric and has bounded support and
bounded derivative. Write
\[
\varphi_K = \int_\R K^2(u)\,\mathrm{d} u\quad
\mbox{and} \quad\psi_K=\frac{1}{2} \int_\R
u^2 K(u) \,\mathrm{d} u.
\]
\end{condition}

\begin{condition}\label{con:regularity}
$\mu,s\in{\cal C}^4([a,b]), \inf_{x\in[a,b]} s(x)>0, f_e\in
{\cal C}^4(\R), f_e(0) > 0, f_e(1)+\break f_e(-1)>0$.
\end{condition}

\subsection{Uniform Bahadur representation for \texorpdfstring{$\hat\mu(x)$}{mu(x)}}\label{sec:hatmu}
Theorem \ref{thmm:m} below provides an asymptotic uniform
Bahadur representation for $\hat\mu(x)$ in (\ref{eq:muest}), and its
proof in Section~\ref{sec:pfthm} relies
on the arguments and results in Section~\ref{sec:coupling}.

\begin{thmm}\label{thmm:m}
Let $\hat\mu(x)$ be as in (\ref{eq:muest}). Assume that
Conditions \ref{con:dc} and \ref{con:reg}--\ref{con:regularity} hold. Further
assume $b_n \to0$ and $(\log N_n)^3/ (N_nb_n) \to0$.
Then
\begin{longlist}[(ii)]
\item[(i)] We have the uniform consistency:
%
%
\begin{equation}
\label{eq:muniform} \sup_{x\in{\cal S}_\epsilon([a,b])} \bigl|\hat\mu(x)-\mu(x)\bigr|
=\mathrm{O}_\mathrm{p} \biggl\{ b^2_n + \frac{(\log N_n)^{3/2}}{(N_nb_n)^{1/2}}
\biggr\}.
\end{equation}
\item[(ii)] Moreover, the Bahadur representation
%
%
\begin{eqnarray}
\label{eq:mbahadur} \hat{\mu}(x)-\mu(x)=\psi_K \rho_\mu(x)
b_n^2 + \frac{(b-a)s(x)}{
f_e(0)} \frac{ Q_{b_n}(x) }{N_nb_n} +
\mathrm{O}_\mathrm{p} ( r_n )
\end{eqnarray}
holds uniformly over $x\in{\cal S}_\epsilon([a,b])$, where
\begin{eqnarray*}
\rho_\mu(x) &=& \mu''(x) - \biggl[
\frac{\mu'(x) f'_e(0)}{f_e(0)} + 2s'(x) \biggr] \frac{\mu'(x)}{s(x)},
\\
Q_{b_n}(x) &=& - \sum^{n}_{i=1}
\sum^{m_i}_{j=1} \bigl\{
\mathbf{1}_{Y_{i,j}\le\mu(x)} - \E[\mathbf{1}_{Y_{i,j}\le\mu
(x)}]\bigr\}
K_{b_n}(x_{i,j}-x),
\\
r_n &=&
b^4_n + \frac{b_n^{1/2} (\log N_n)^{3/2}}{N_n^{1/2}} + \frac{(\log N_n)^{9/4}}{(N_nb_n)^{3/4}}.
\end{eqnarray*}
\end{longlist}
\end{thmm}

In the Bahadur representation (\ref{eq:mbahadur}), $\psi_K
\rho_\mu(x) b_n^2$ is the bias term, $Q_{b_n}(x)$ determines
the asymptotic distribution of $\hat{\mu}(x)-\mu(x)$, and $r_n$
is the negligible error term. Such a Bahadur representation
provides a powerful tool in studying the asymptotic behavior of
$\hat\mu(x)$. Based on Theorem \ref{thmm:m}, we obtain a Central
Limit theorem (CLT) for $\hat\mu$ in Corollary \ref{cor:m}.
Clearly, the variance of $Q_{b_n}(x)$ is a linear
combination of $K_{b_n}(x_{i,j_1}-x)K_{b_n}(x_{i,j_2}-x)$. The
following regularity condition is needed to ensure the
negligibility of the cross-term
$K_{b_n}(x_{i,j_1}-x)K_{b_n}(x_{i,j_2}-x)$ for $j_1\ne j_2$.

\begin{condition}\label{con:clt}
Assume that, for all given $x\in{\cal S}_\epsilon([a,b])$ and
$k_n=\mathrm{O}(\log N_n)$, there exits $\iota_n$ such that $k_n
\iota_n\to0$ and
%
%
\begin{equation}
\label{eq:conclt} \sum_{(i,j_1,j_2)\in{\cal I}} K_{b_n}(x_{i,j_1}-x)K_{b_n}(x_{i,j_2}-x)
= \mathrm{O}\bigl[ \min(h,M_n)nb_n k_n
\iota_n\bigr],\qquad  M_n=\max_{1\le i\le n}
m_i\quad
\end{equation}
for all $h\ge(k_n\vee a)$, where ${\cal I}=\{(i, j_1,j_2)\dvt 1\le i\le n, a\le j_1 < j_2\le\min(a+h-1,m_i), |j_1-j_2|\le
k_n\}$. Further assume that $\max_j
\sum^n_{i=1}K^r_{b_n}(x_{i,j}-x)=\mathrm{O}(nb_n),r=2,4$.
\end{condition}

Condition \ref{con:clt} is very mild. Intuitively, we
consider $x_{i,j}, j\in\Z$, being random locations,
then
$\E[K_{b_n}(x_{i,j_1}-x)K_{b_n}(x_{i,j_2}-x)]=\mathrm{O}(b^2_n)$ for
$j_1\ne j_2$.
Thus, under the mild condition
$b_n \log N_n\to0$, (\ref{eq:conclt}) holds with $\iota_n=b_n$.

\begin{corol}\label{cor:m}
Let the conditions in Theorem \ref{thmm:m} be fulfilled and
Condition \ref{con:clt} hold. Further assume that $(\log
N_n)^9/(N_nb_n) + N_n b_n^9 \to0$ and $nM_n=\mathrm{O}(N_n),
nb_n\to\infty, \log N_n=\mathrm{O}(\sqrt{M_n})$, where $M_n$ is defined
as in (\ref{eq:conclt}). Then, for any $x\in{\cal
S}_\epsilon([a,b])$, we have
%
%
\begin{equation}
\label{cor:mconclt} (N_nb_n)^{1/2} \bigl[ \hat{
\mu}(x)-\mu(x) - \psi_K \rho_\mu(x) b_n^2
\bigr] \Rightarrow N \biggl(0, \frac{ \varphi_K (b-a)s^2(x) }{4
f^2_e(0) } \biggr).
\end{equation}
\end{corol}

The proof of Corollary \ref{cor:m}, given in Section~\ref{sec:clt},
uses the coupling argument in Section~\ref{sec:coupling}.
The condition $nM_n=\mathrm{O}(N_n)$ is in line with the classical
CLT Lindeberg condition that none of the subjects dominates the
others. If $b_n$ is of the order $N_n^{-\beta}$, then the bandwidth
condition in Corollary \ref{cor:m} holds if $ \beta\in
(1/9,1)$. By Corollary
\ref{cor:m}, the optimal bandwidth minimizing the asymptotic mean
squared error is
%
%
\begin{equation}
\label{eq:bndwth} b_n = \biggl[ \frac{ \varphi_K (b-a) s^2(x) } { 4 \psi_K^2
\rho^2_\mu(x) f_e^2 (0)}
\biggr]^{1/5} N^{-1/5}_n.
\end{equation}
For this optimal bandwidth, the bias term is of order
$\mathrm{O}(N^{-2/5}_n)$ and contains the derivatives $s',\mu',\mu''$
and $f'_e$ that can be difficult to estimate. Based on the Bahadur
representation
(\ref{eq:mbahadur}), we can correct the bias term $\psi_K
\rho_\mu(x) b_n^2$ via the jackknife estimator
$\tilde{\mu}(x)$ in (\ref{eq:mbias}). Then the bias term for
$\tilde{\mu}(x)$ becomes $2\psi_K \rho_\mu(x) b_n^2
-\psi_K \rho_\mu(x) (\sqrt{2}b_n)^2=0$. By
(\ref{eq:mbahadur}), following the proof of Corollary
\ref{cor:m}, we have
%
%
\begin{equation}
\label{eq:mjackclt} (N_nb_n)^{1/2} \bigl[
\tilde{\mu}(x)-\mu(x) \bigr] \Rightarrow N \biggl(0, \frac{ \varphi_{K^*} (b-a) s^2(x)} {
4f^2_e(0) } \biggr),
\end{equation}
where $K^*(u) = 2K(u)-2^{-1/2}K(u/\sqrt{2})$.

\subsection{Uniform Bahadur representation for \texorpdfstring{$\hat{s}(x)$}{s(x)}}\label{sec:vbahadur}
Theorem \ref{thmm:v} below provides a uniform
Bahadur representation for $\hat{s}(x)$ in (\ref{eq:sest}).

\begin{thmm}\label{thmm:v}
Let $\hat{s}(x)$ be as in (\ref{eq:sest}).
Assume that the conditions in Theorem \ref{thmm:m} hold. Further
assume ${h_n} + (\log N_n)^3 /(N_n h_n) \to0$.
Then
\begin{longlist}[(ii)]
\item[(i)] We have the uniform consistency:
%
%
\begin{eqnarray}
\label{eq:muniform} \sup_{x\in{\cal S}_\epsilon([a,b])} \bigl|\hat{s}(x)-s(x)\bigr|
=\mathrm{O}_\mathrm{p} \biggl\{ b^2_n +
h^2_n + \frac{(\log N_n)^{3/2}}{(N_nb_n)^{1/2}} + \frac{(\log N_n)^{3/2}}{(N_nh_n)^{1/2}} \biggr\}.
\end{eqnarray}
\item[(ii)] Moreover, the Bahadur representation
%
%
\begin{eqnarray}
\label{eq:vbahadur} \hat s(x) -s(x) = \psi_K \rho_s(x)
{h_n^2} + (b-a)s(x) \biggl[ \frac{ W_{h_n} (x) }{N_nh_n\kappa_+ } -
\frac{\kappa T_{b_n}(x) } { N_nb_n f_{e}(0) } \biggr] + \mathrm{O}_\mathrm{p}(\tilde{r}_n),
\end{eqnarray}
holds uniformly over $x\in{\cal S}_\epsilon([a,b])$, where
$\kappa_{+} = f_e(-1) + f_e(1),
\kappa=[f_e(1)-f_e(-1)]/\kappa_+$, $Q_{b_n} (x)$ is defined as in
Theorem \ref{thmm:m},
\begin{eqnarray*}
T_{b_n}(x) &=& 2Q_{b_n} (x)- 2^{-1/2}
Q_{\sqrt{2}b_n}(x),
\\
\rho_s(x) &=& s''(x)
- \frac{2s'(x)^2}{s(x)}+ \kappa \biggl[ \mu''(x) -
\frac{2\mu'(x)s'(x)}{s(x)} \biggr]
\\
&&{} - \frac{
f'_{e}(1) [ s'(x)+\mu'(x) ]^2 -f'_{e}(-1) [ s'(x)-\mu'(x)
]^2}{\kappa_+ s(x)},
\\
W_{h_n} (x) &
= & - \sum^{n}_{i=1} \sum
^{m_i}_{j=1} \bigl\{ \mathbf{1}_{|Y_{i,j}-\mu(x)|\le s(x)} - \E[
\mathbf{1}_{|Y_{i,j}-\mu(x)|\le s(x)}] \bigr\} K_{h_n}(x_{i,j}-x),
\\
\tilde{r}_n &=& b^4_n + h^4_n
+ \frac{h_n^{1/2} (\log N_n)^{3/2}}{N_n^{1/2}} + \frac{(\log N_n)^{9/4}}{(N_nh_n)^{3/4}}
\\
&&{} + \frac{(\log N_n)^{9/4}}{N_n^{3/4} b_n^{1/4} h_n^{1/2} } +
\frac{b_n (\log N_n)^{3/2}}{(N_n h_n)^{1/2}}.
\end{eqnarray*}
\end{longlist}
\end{thmm}

As in Corollary \ref{cor:m}, we can use the Bahadur
representation (\ref{eq:vbahadur}) to obtain a CLT for $\hat
s(x) -s(x)$. However, the convergence rate depends on the
ratio $h_n/b_n$. If $h_n/b_n\to\infty$, then the term
$T_{b_n}(x)/(N_nb_n)$ dominates and we have
$(N_nb_n)^{1/2}$-convergence; if $h_n/b_n\to0$, then the term
$W_{h_n}(x)/(N_nh_n)$ dominates and we have
$(N_nh_n)^{1/2}$-convergence; if $h_n/b_n \to c$ for a constant
$c\in(0,\infty)$, then both terms contribute.

\begin{corol}\label{cor:v}
Let the conditions in Theorem \ref{thmm:v} be fulfilled and
Condition \ref{con:clt} and its counterpart version with $b_n$
being replaced by $h_n$ hold. Further assume that
\[
N_n (b_n\vee h_n)^9 +
\frac{(\log N_n)^9}{N_n(b_n\wedge h_n)} \to0,
\]
and $nM_n=\mathrm{O}(N_n), n(b_n\wedge h_n) \to\infty, \log
N_n=\mathrm{O}(\sqrt{M_n})$, where $M_n$ is defined as in
(\ref{eq:conclt}). Recall $K^*(u) =
2K(u)-2^{-1/2}K(u/\sqrt{2})$ in (\ref{eq:mjackclt}) and
$\kappa,\kappa_+$ in Theorem \ref{thmm:v}. Let $x\in{\cal
S}_\epsilon([a,b])$ be a fixed point. Suppose
$h_n/b_n\to c$.
\begin{longlist}[(iii)]
\item[(i)] If $\kappa\ne0$ and $c=\infty$, then
\[
(N_nb_n)^{1/2} \bigl[ \hat{s}(x)-s(x) -
\psi_K \rho_s(x) h_n^2 \bigr]
\Rightarrow N \biggl(0, \frac{ \varphi_{K^*} \kappa^2 (b-a) s^2(x) }{4
f^2_e(0) } \biggr).
\]

\item[(ii)] If $\kappa\ne0$ and $c\in[0,\infty)$, then
%
%
\begin{equation}
\label{cor:vcltii} (N_nh_n)^{1/2} \bigl[
\hat{s}(x)-s(x) - \psi_K \rho_s(x) h_n^2
\bigr] \Rightarrow N \bigl(0, \sigma_c^2\bigr),
\end{equation}
where
\[
\sigma^2_c = \frac{(b-a)s^2(x)}{4} \biggl\{
\frac{ \varphi_{K} } {\kappa_{+}^2 } + \frac{c^2 \kappa^2 \varphi_{K^*} }{f^2_e(0)} -\frac{ 2 c \kappa[ 1-4F_e(-1) ] } {
\kappa_{+} f_e(0)} \int
_{\R} K(u) K^*(c u) \,\mathrm{d}u \biggr\}.
\]

\item[(iii)] If $\kappa= 0$, then for all $c \in[0,\infty]$,
(\ref{cor:vcltii}) holds with $\sigma_c^2 =
\varphi_{K} (b-a) s^2(x) / ( 4 \kappa_{+}^2 )$.

\end{longlist}
\end{corol}

One can similarly establish CLT results for $\tilde{s}(x)$ in
(\ref{eq:vbias}). We omit the details.

\section{An illustration using real data}\label{sec:app}

\subsection{Bandwidth selection}
For least-squares based estimation of longitudinal data,
Rice and Silverman \cite{ric1991} suggested the subject-based
cross-validation method. The basic idea is to use all but one subject
to do model fitting, validate the
fitted model using the left-out subject, and finally choose the
optimal bandwidth by minimizing the overall prediction error:
%
%
\begin{equation}
\label{eq:SJCV} b^*_\mathrm{LS}=\argmin_{b} \sum
_{i=1}^n \sum_{j=1}^{m_i}
\bigl\{ Y_{i,j}-\tilde\mu^{(-i)}(x_{i,j})\bigr
\}^2,
\end{equation}
where $\tilde\mu^{(-i)}(x)$ represents the estimator of $\mu(x)$
based on data from all but $i$th subject. As in Wei, Zhao and Lin \cite
{wei2012}, we replace the square loss by
absolute deviation:
%
%
\begin{equation}
\label{eq:SJCV2} b^*_\mathrm{LAD}=\argmin_{b} \sum
_{i=1}^n \sum_{j=1}^{m_i}
\bigl|Y_{i,j}-\tilde\mu^{(-i)}(x_{i,j})\bigr|.
\end{equation}

\subsection{An illustration using progesterone data}\label{sec:app2}
\begin{figure}

\includegraphics{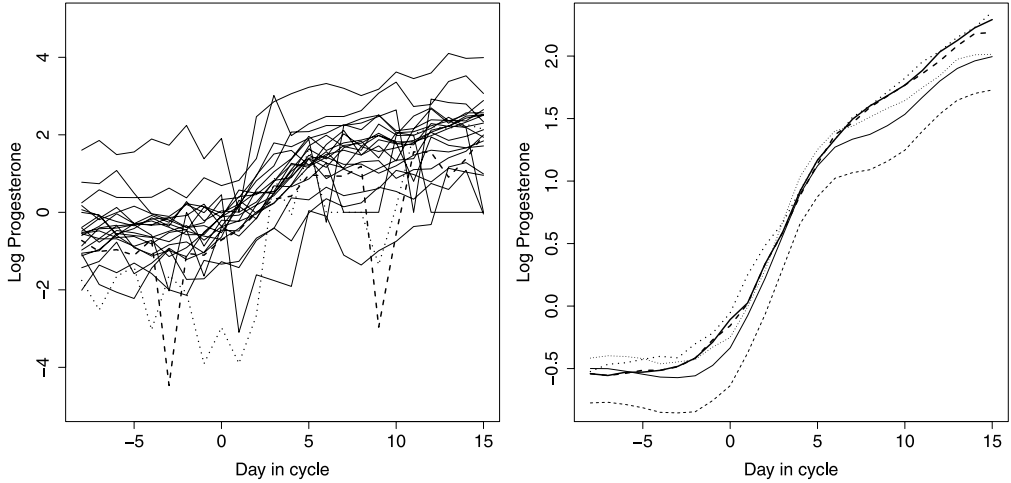}

\caption{Left: Trajectories of the measurements from 22 conceptive
women. Right: Estimates of $\mu(\cdot)$
using both the original data and perturbed data.
Thin solid, dotted, and dashed curves are the least-squares estimates
of $\mu(\cdot)$ based on the original data,
perturbation scenario I (remove subjects 13 and 14), and perturbation
scenario II (shift subjects 13 and 14 down by three units), respectively.
Similarly, thick solid, dotted, and dashed curves are
least-absolute-deviation estimates.}\label{fig1}
\end{figure}

Urinary metabolite progesterone levels are measured daily,
around the ovulation day,
over 22 conceptive and 69 nonconceptive women's menstrual cycles so
that each curve has
about 24 design points; see the left panel of Figure~\ref{fig1} for a plot of
the trajectories of the 22 conceptive women.
Previous studies based on least-squares (LS) methods include
Brumback and Rice \cite{bru1998}, Fan and Zhang \cite{fan2000}, and Wu
and Zhang \cite{wu2002}.
Here we reanalyze the conceptive group using our
least-absolute-deviation (LAD) estimates.

From the left plot in Figure~\ref{fig1},
subject 14 (dashed curve) has two sharp drops in progesterone levels at
days $-3$ and 9. Similarly,
subject 13 (dotted curve) has unusually low levels on days $-1,0,1$.
While such sharp drops or ``outliers'' may be caused by incorrect
measurements or other
unknown reasons, we investigate the impact of such ``outliers'' on the
LS and LAD estimates.
In the right plot of
Figure~\ref{fig1}, the thick solid and thin solid curves are the LAD and LS
estimates of $\mu(\cdot)$. The two estimates are
reasonably close except during the periods $[-4,1]$ and $[8,15]$.
Notice that the latter periods contain the ``outliers'' from subjects
13, 14.

To understand the impact of such possible ``outliers'', we consider two
scenarios of perturbing the data below.
\begin{longlist}[(ii)]
\item[(i)] Scenario I: remove subjects 13 and 14 and estimate
$\mu(\cdot)$ using the remaining subjects. The thick dotted and thin
dotted curves are the corresponding LAD and LS estimates. Clearly, the
discrepancy is
largely diminished.

\item[(ii)] Scenario II: make the two outlier subjects 13 and 14 even
more extreme by shifting their curves three units down. We see that the
discrepancy
between the LAD (thick dashed) and LS (thin dashed) estimates becomes
even more remarkable.
\end{longlist}
Compared with the estimate based on the original data,
the LS estimates under the two perturbation scenarios differ significantly.
By contrast, the LAD estimates under the three cases are similar,
indicating the robustness
in the presence of outliers. We conclude that, for the progesterone
data with several possible outliers,
the proposed LAD estimate offers an attractive alternative over the
well-studied LS estimates.
In practice, we recommend the LAD estimate if the data has suspicious,
unusual observations or extreme outliers.

\section{Conclusion and extension to spatial setting}\label{sec:dis}

This paper studies robust estimations of the location and scale
functions in
a nonparametric regression model with serially dependent data from
multiple subjects.
Under a general error dependence structure that allows for many linear
and nonlinear processes,
we study uniform Bahadur representations
and asymptotic normality for
least-absolute-deviation estimations of a location-scale longitudinal model.
In the large literature on nonparametric estimation of
longitudinal models, most existing works use least-squares based
methods, which
are sensitive to extreme observations and may perform poorly
in such circumstances.
Despite the popularity of quantile regression methods in linear models
and nonparametric
regression models, little research has been done in quantile regression
based estimations for nonparametric
longitudinal models, partly due to difficulties in dealing with the dependence.
Therefore, our work provides a solid theoretical foundation for
quantile regression estimations in longitudinal models.

The study of asymptotic Bahadur representations is a difficult area and
has focused
mainly on the i.i.d. setting or stationary time series setting. For longitudinal
data, deriving Bahadur representations is more challenging due to
the nonstationarity and dependence.
To obtain our Bahadur representations, we develop substantial theory
for kernel weighted empirical processes via a coupling argument.

The proposed error dependence structure and coupling argument provide a
flexible and powerful framework for
asymptotics from dependent data, such as time series data, longitudinal
data and spatial data, whereas similar problems have
been previously studied mainly for either independent data or
stationary time series.
In (\ref{eq:error}), $e_j$
depends on the innovations or shocks $\varepsilon_j,\varepsilon_{j\pm
1},\ldots,$ indexed by integers on a line.
A natural extension is the function of innovations indexed by bivariate
integers on a square:
\[
e_j = G(\varepsilon_{j,j}, \varepsilon_{j,j\pm1},
\varepsilon_{j\pm
1,j}, \varepsilon_{j\pm1, j\pm1}, \ldots),\qquad j\in\Z.
\]
The coupling argument still holds by replacing the innovations
$\varepsilon_{j\pm r,j\pm s}, r,s \ge k+1$, outside the $k$ nearest squares
with i.i.d. copies. As in Condition \ref{con:dc}, we can assume that
the impact of perturbing the distant innovations decays exponentially
fast (or polynomially fast with slight modifications of the proof).
More generally, the coupling argument holds for function of
innovations indexed by multivariate spatial lattice, and such setting
may be useful in studying asymptotics for spatial data.

\section{Technical proofs}\label{sec:proof}
Throughout $c,c_1,c_2,\ldots,$ are generic constants. First, we give an
inequality for the indicator function.
Let $Z,Z'$ be two random variables and $y\in\R$. For $\alpha>0$, we have
\[
\mathbf{1}_{Z\le y< Z'} = \mathbf{1}_{Z\le y < Z', |Z-Z'|\ge\alpha} +
\mathbf{1}_{Z\le y < Z', |Z-Z'| < \alpha} \le\mathbf{1}_{|Z-Z'|\ge
\alpha} + \mathbf{1}_{y < Z' < y + \alpha},
\]
similarly, $\mathbf{1}_{Z' \le y < Z}
\le\mathbf{1}_{|Z-Z'|\ge
\alpha} + \mathbf{1}_{y -\alpha< Z' \le y }$.
Therefore,
%
%
\begin{equation}
\label{eq:zineq} |\mathbf{1}_{Z\le y} - \mathbf{1}_{Z'\le y}| =
\mathbf{1}_{Z\le y< Z'}+\mathbf{1}_{Z' \le y <
Z} \le2 \mathbf{1}_{|Z-Z'| \ge\alpha}
+ \mathbf{1}_{ y-\alpha< Z' < y+\alpha}.
\end{equation}

\subsection{Proof of Propositions \texorpdfstring{\protect\ref{pro:h}--\protect\ref{pro:e}}{2.1--2.2}}
\begin{pf*}{Proof of Proposition \ref{pro:h}}
Let $q^*=q/(\varsigma+\upsilon),
p_1=\upsilon/\varsigma+1$, and $p_2=\varsigma/\upsilon+1$ so that
$\varsigma q^* p_1=q, \upsilon q^* p_2 = q$, and $1/p_1+1/p_2=1$.
For convenience, write $e'_0=e_0(k)$. By assumption, $\|e_0-e'_0\|
_q=\mathrm{O}(\rho^k)$.
By (\ref{eq:Hfunc}) and the H\"older inequality $\E|Z_1Z_2|\le
\|Z_1\|_{p_1} \|Z_2\|_{p_2}$,
\begin{eqnarray*}
\bigl\|h\bigl(e'_0\bigr)-h(e_0)
\bigr\|^{q^*}_{q^*} & \le& \mathrm{O}(1) \E\bigl[ \bigl|e'_0-e_0\bigr|^{\varsigma
q^*}
\bigl(1+|e_0|+\bigr|e'_0\bigr|\bigr)^{\upsilon q^*}
\bigr]
\cr
&\le& \mathrm{O}(1) \bigl\{\E\bigl[\bigl|e_0-e'_0\bigr|^{\varsigma q^* \cdot p_1}
\bigr]\bigr\}^{1/p_1} \bigl\{ \E\bigl[\bigl(1+|e_0|+\bigl|e'_0\bigr|
\bigr)^{\upsilon q^* \cdot p_2}\bigr]\bigr\}^{1/p_2}
\cr
&=& \mathrm{O}(1)
\bigl\|e_0-e'_0\bigr\|_{q}^{q/p_1}
\|e_0\|_{q}^{q/p_2}=\mathrm{O}\bigl(\rho^{k q/p_1}
\bigr).
\end{eqnarray*}
The above expression gives $\|h(e'_0)-h(e_0)\|_{q^*}\le \mathrm{O}(1) [\rho
^{q/(p_1 q^*)}]^k=\mathrm{O}(\rho^{k\varsigma})$.
\end{pf*}

\begin{pf*}{Proof of Proposition \ref{pro:e}}
Let $\alpha=\rho^{kq/(1+q)}$. By (\ref{eq:zineq}) and the triangle inequality,
\begin{eqnarray*}
\|\mathbf{1}_{e_0\le x}-\mathbf{1}_{e_0(k)\le x}\|_q &\le& 2
\|\mathbf{1}_{|e_0-e_0(k)| \ge\alpha}\|_q + \|\mathbf {1}_{x-\alpha\le e_0 \le x+\alpha}
\|_q
\cr
&=& 2 \bigl[ \p\bigl\{\bigl|e_0-e_0(k)\bigr|
\ge\alpha\bigr\} \bigr]^{1/q} + \bigl[\p\{x-\alpha\le e_0
\le x+\alpha\}\bigr]^{1/q}.
\end{eqnarray*}
By the Markov inequality, $\p\{|e_0-e_0(k)| \ge\alpha\}\le\E
[|e_0-e_0(k)|^q]/\alpha^q=\mathrm{O}(\rho^{kq}/\alpha^q)$.
Since $e_0$ has a bounded density, $\p\{x-\alpha\le e_0 \le x+\alpha\}
=\mathrm{O}(\alpha)$. The result then follows.
\end{pf*}

\subsection{Proof of Theorems \texorpdfstring{\protect\ref{thmm:cp}--\protect\ref{pro:osc}}{3.1--3.3}} \label{sec:approsc}

\begin{pf*}{Proof of Theorem \ref{thmm:cp}}
For $r=1,2,\ldots,2k_n+1$, let
%
%
\begin{equation}
\label{eq:is} {\cal I}_r=\bigl\{(i,j)\dvt 1\le i\le n, 1\le j\le\bigl
\lfloor(m_i-r)/(2k_n+1) \bigr\rfloor+ 1 \bigr\}.
\end{equation}
Using the identity $\sum^m_{j=1} a_j = \sum^{k}_{r=1} \sum^{\lfloor
(m-r)/k \rfloor+ 1 }_{j=1} a_{(j-1)k+r}$ for all $k,m\in\N, a_1,\ldots
,a_m\in\R$, we can rewrite $H_n(v)$ as
%
%
\begin{equation}
\label{eq:hnd} H_n(v)=\sum^{2k_n+1}_{r=1}
\sum_{(i,j)\in{\cal I}_r} g_{i,(j-1)(2k_n+1)+r}(\tilde{Y}_{i,(j-1)(2k_n+1)+r},v)
:= \sum^{2k_n+1}_{r=1} H_n(v,r).
\end{equation}
Now we consider $H_n(v,r)$. By the discussion in Section~\ref{sec:coupling},
the
summands in $H_n(v,r)$ are independent. By (\ref{eq:thmcpcon}),
%
%
\begin{eqnarray}
\label{eq:varh} \operatorname{Var}\bigl[H_n(v,r)\bigr] &=&\sum
_{(i,j)\in{\cal I}_r} \E\bigl[ g^2_{i,(j-1)(2k_n+1)+r}(
\tilde{Y}_{i,(j-1)(2k_n+1)+r},v) \bigr]
\nonumber
\\[-8pt]
\\[-8pt]
\nonumber
& \le& \sum^{n}_{i=1}
\sum^{m_i}_{j=1} \E\bigl[
g^2_{i,j}(\tilde{Y}_{i,j},v) \bigr] \le
\chi_n,
\end{eqnarray}
uniformly over $v,r$.

(i) Consider the case $\chi_n=\mathrm{O}(1)$. Recall the condition
$|g_{i,j}(y,v)|\le c$. By
Berstein's exponential inequality (Bennett \cite{ben1962}) for bounded
and independent random variables, for any given $c_1>0$, when $N_n$ is
sufficiently large,
%
%
\begin{eqnarray}
\label{eq:berstein} \p\bigl\{ \bigl|H_n(v,r)\bigr|\ge c_1 \log
N_n \bigr\} \le 2 \exp \biggl\{-\frac{(c_1 \log N_n)^2}{2\operatorname{Var}[\Lambda
_n(r,h)] + c c_1 \log N_n } \biggr\} \le2
N_n^{-c_1/(3c)},
\end{eqnarray}
uniformly over $r$ and $h$. Here the second inequality follows from
$\operatorname{Var}[H_n(v,r)]\le\chi_n=\mathrm{O}(1)\le cc_1 \log N_n$
for large enough $N_n$.
Thus,
\begin{eqnarray*}
\p \Bigl\{ \max_{v\in\mathcal{V}_n, 1\le r\le2k_n+1} \bigl|H_n(v,r)\bigr| \ge
c_1 \log N_n \Bigr\} &\le&\sum
_{v\in\mathcal{V}_n, 1\le r\le2k_n+1} \p\bigl\{\bigl|H_n(v,r)\bigr|\ge c_1
\log N_n \bigr\}
\\
&\le& 2 |\mathcal{V}_n|
k_n N_n^{-c_1/(3c)}.
\end{eqnarray*}
By the assumption that both $|\mathcal{V}_n|$ and $k_n$ grow no faster
than a polynomial of $N_n$, we can make the above probability go to
zero by choosing a large enough $c_1$.
Therefore, $\max_{v\in\mathcal{V}_n, 1\le r\le2k_n+1}
|H_n(v,r)|=\mathrm{O}_\mathrm{p}(\log N_n)$. By (\ref{eq:hnd}), the desired
result follows
from
\[
\max_{v\in\mathcal{V}_n} \bigl|H_n(v)\bigr|\le(2k_n+1) \max
_{v\in\mathcal{V}_n,
1\le r\le2k_n+1}\bigl |H_n(v,r)\bigr|.
\]

(ii) Consider the case $\sup_n \log N_n/\chi_n<\infty$. As in (\ref
{eq:berstein}),
\begin{eqnarray*}
\p\bigl\{ \bigl|H_n(v,r)\bigr|\ge c_1 \sqrt{\chi_n
\log N_n} \bigr\} \le 2 \exp \biggl\{-\frac{(c_1 \sqrt{\chi_n \log N_n})^2}{2\chi_n + c
c_1 \sqrt{\chi_n \log N_n} } \biggr\} =\mathrm{O}
\bigl[N_n^{-c_1^2/(2+cc_1c_2)}\bigr],
\end{eqnarray*}
uniformly over $r$ and $h$, where $c_2=\sup_n [ \log N_n/\chi_n
]^{1/2}<\infty$. The rest of the proof follows from the same argument
as in the case (i) by
choosing a sufficiently large $c_1$.
\end{pf*}

\begin{pf*}{Proof of Theorem \ref{lem:appr}}
Let $\alpha=1/N_n$. Since $\varpi_{i,j}(x)\le c$, applying (\ref
{eq:zineq}), we obtain
%
%
\begin{eqnarray}
\label{eq:omegan} \bigl|F_n(x,y)-\tilde{F}_n(x,y)\bigr| & \le& \sum
^{n}_{i=1} \sum
^{m_i}_{j=1} \varpi_{i,j}(x) |
\mathbf{1}_{Y_{i,j}\le y} - \mathbf{1}_{\tilde{Y}_{i,j}\le y} |
\nonumber\\
& \le& 2c \Biggl[
\sum^{n}_{i=1} \sum
^{m_i}_{j=1} \mathbf {1}_{|Y_{i,j}-\tilde{Y}_{i,j}|\ge\alpha} + \sum
^{n}_{i=1} \sum
^{m_i}_{j=1} \mathbf{1}_{y-\alpha<\tilde
{Y}_{i,j}<y+\alpha} \Biggr]
\\
&:=& 2c \bigl[ \Omega_n + \Lambda_n(y) \bigr].\nonumber
\end{eqnarray}
Notice that, $|Y_{i,j}-\tilde{Y}_{i,j}| =\mathrm{O}(1) |e_{i,j}-e_{i,j}(k_n)|$.
By (\ref{eq:gmc}) and
the Markov inequality,
\[
\E(\mathbf{1}_{|Y_{i,j}-\tilde{Y}_{i,j}|\ge\alpha}) \le \frac{ \|Y_{i,j}-\tilde{Y}_{i,j}\|^q_q}{\alpha^q} = \mathrm{O}(1) \frac{\|e_{i,j}-e_{i,j}(k_n)\|^q_q}{\alpha^q}
=\mathrm{O}\bigl(N_n^q \rho^{qk_n}\bigr).
\]
Thus, $\Omega_n=\mathrm{O}_\mathrm{p}(N_n^{1+q} \rho^{qk_n})=\mathrm{O}_\mathrm
{p}[N_n^{1+q} \rho^{q\lambda\log(N_n)}]=\mathrm{o}_\mathrm{p}(1)$ for
$\lambda>(q+1)/[q\log(1/\rho)]$.

For $\Lambda_n(y)$ over $y\in\R$, consider two cases: $|y|>N_n^{1/q}$
and $|y|\le N_n^{1/q}$.
For $|y|>N_n^{1/q}$, since $\alpha=1/N_n\to0$, $\mu(x_{i,j})$ and
$s(x_{i,j})$ are bounded,
$\{y-\alpha<\tilde{Y}_{i,j}<y+\alpha\}\subset\{|e_{i,j}(k_n)|\ge c_1
N_n^{1/q}\}$ for some constant $c_1>0$.
Therefore, by $e_{i,j}(k_n)\in{\cal L}^q$ and the Markov inequality,
%
%
\begin{eqnarray}
\label{eq:pthm1a} \E \Bigl[ \sup_{|y|>N_n^{1/q}} \Lambda_n(y)
\Bigr]&\le&\E \Biggl[ \sum^{n}_{i=1} \sum
^{m_i}_{j=1} \mathbf {1}_{|e_{i,j}(k_n)|>c_1 N_n^{1/q}}
\Biggr]
\nonumber
\\[-8pt]
\\[-8pt]
\nonumber
 &\le&\sum^{n}_{i=1} \sum
^{m_i}_{j=1} \frac{\|e_{i,j}(k_n)\|^q_q}{(c_1
N_n^{1/q})^q}=\mathrm{O}(1).
\end{eqnarray}
We conclude that $\sup_{|y|>N_n^{1/q}} \Lambda_n(y)=\mathrm{O}_\mathrm{p}(1)$.

In what follows, we use a chain argument to prove $\sup_{y\in
[-N_n^{1/q},N_n^{1/q}]} \Lambda_n(y)=\mathrm{O}_\mathrm{p}[(\log n)^2]$.
Without loss of generality, consider $y\in[0,N_n^{1/q}]$.
Write $\ell_n=\lfloor N_n^{1+1/q} \rfloor$ and let $\mathcal{V}_n=\{
y_v=v N_n^{1/q}/\ell_n, v=0,1,\ldots,\ell_n\}$ be the set of
$\ell_n+1$ grid
points uniformly spaced over $[0,N_n^{1/q}]$.
Partition $[0,N_n^{1/q}]$ into intervals $I_v=[y_{v-1},y_v], v=1,\ldots
,\ell_n$. For any
$y\in I_v$, we have $\mathbf{1}_{y-\alpha
<\tilde{Y}_{i,j}<y+\alpha} \le\mathbf{1}_{y_{v-1}-\alpha
<\tilde{Y}_{i,j}<y_v+\alpha}$. Since $s(x_{i,j})$ is bounded away from
zero, $\sup_{u}f_e(u)<\infty$, and $|y_v-y_{v-1}|=\mathrm{O}(1/N_n)$, we have
$\E({\bf
1}_{y_{v-1}-\alpha<\tilde{Y}_{i,j}<y_v+\alpha})\le c_2/N_n$
uniformly for some constant $c_2<\infty$. Consequently, for
any $y\in I_v$, we have
\begin{eqnarray*}
\Lambda_n(y) \le \sum^{n}_{i=1}
\sum^{m_i}_{j=1} \bigl[ \bigl\{
\mathbf{1}_{y_{v-1}-\alpha<\tilde{Y}_{i,j}<y_v+\alpha} - \E(\mathbf{1}_{y_{v-1}-\alpha<\tilde{Y}_{i,j}<y_v+\alpha}) \bigr\} +
c_2/N_n \bigr] = \Lambda^*_n(v) +
c_2.
\end{eqnarray*}
We apply Theorem \ref{thmm:cp} to $\Lambda^*_n(v)$. For $\chi_n$ in
(\ref{eq:thmcpcon}), using $\E({\bf
1}_{y_{v-1}-\alpha<\tilde{Y}_{i,j}<y_v+\alpha})\le c_2/N_n$, we have
$\chi_n=\mathrm{O}(1)$
and thus $\max_{v\in\mathcal{V}_n}|\Lambda_n^*(v)|=\mathrm{O}_\mathrm{p}[(\log
N_n)^2]$, completing the proof.
\end{pf*}

\begin{pf*}{Proof of Theorem \ref{pro:osc}}
Recall the coupling process $\tilde{D}_n(\delta,x,y)$ in (\ref
{eq:empdiff2}). Under the assumption $\sup_n \log N_n/(\delta_n\phi
_n)<\infty$,
$(\log N_n)^2=\mathrm{O}\{ [\delta_n\phi_n (\log N_n)^3]^{1/2}\}$. Thus, by
(\ref{eq:dd2}),
it suffices to show $\sup_{|\delta|\le\delta_n, x,y\in\R} |\tilde
{D}_n(\delta,x,y)|=\mathrm{O}_\mathrm{p}\{ [\delta_n\phi_n (\log N_n)^3]^{1/2}\}$.

Without loss of generality, assume $\delta\in[0,\delta_n]$.
Recall $\tilde{Y}_{i,j}$ in (\ref{eq:Fn2}).
Rewrite
\begin{eqnarray*}
\tilde{D}_n(\delta,x,y) = \sum^n_{i=1}
\sum^{m_i}_{j=1}\varpi_{i,j}(x)
\bigl\{ \tilde\xi_{i,j}(\delta ,y) - \E\bigl[ \tilde\xi_{i,j}(
\delta,y) \bigr]\bigr\},\qquad \tilde\xi_{i,j}(\delta,y) =
\mathbf{1}_{y<\tilde{Y}_{i,j}\le
y+\delta}.
\end{eqnarray*}
As in the proof of Theorem \ref{lem:appr}, consider
$|y|>N_n^{1/q}$ and $|y|\le N_n^{1/q}$.

For $|y|>N_n^{1/q}$, since $\mu(x_{i,j})$ and $s(x_{i,j})$ are bounded
and $|\delta|\le\delta_n\to0$,
$\{y<\tilde{Y}_{i,j}\le
y+\delta\} \subset\{|e_{i,j}(k_n)|\ge c_1 N_n^{1/q}\}$ for some
$c_1>0$. Therefore, by the boundedness of
$\varpi_{i,j}(\cdot)$, the same argument in (\ref{eq:pthm1a}) shows
$\tilde{D}_n(\delta,x,y)=\mathrm{O}_\mathrm{p}(1)$ uniformly over $x\in\R, |y|>
N_n^{1/q},\break  |\delta|\le\delta_n$.

Next, we consider $|y|\le N_n^{1/q}$.
Since $\varpi_{i,j}(x)$ vanishes for $x$ outside a bounded interval,
without loss of generality we only
consider $x\in[0,b]$ for some $b>0$, $y\in[0,N_n^{1/q}]$, and $\delta
\in[0,\delta_n]$.
As in the proof of Theorem \ref{lem:appr}, we use the chain argument. Let
$\ell_n=\lfloor N_n^{1/q}/\delta_n + N_n\tau_n +N_n^{1+1/q}\rfloor$,
and
\begin{eqnarray*}
\mathcal{V}_n= \biggl\{(x_{v_1},y_{v_2},t_{v_3}),
x_{v_1}=\frac{v_1
b}{\ell_n}, y_{v_2}=\frac{v_2 N_n^{1/q}}{\ell_n},t_{v_3}=
\frac{v_3 \delta_n}{\ell_n}, v_1,v_2,v_3=0,1,\ldots,
\ell_n \biggr\}
\end{eqnarray*}
be uniformly spaced grid points.
Partition $[0,b]\times[0,N_n^{1/q}]\times
[0,\delta_n]$ into intervals $I_{v_1,v_2,v_3}=[x_{v_1-1},x_{v_1}]\times
[y_{v_2-1},y_{v_2}]\times[t_{v_3-1},t_{v_3}], v_1,v_2,v_3=1,\ldots,
\ell_n$.
Let
\begin{eqnarray*}
\underline{\xi}_{i,j}(v_2,v_3) =
\mathbf{1}_{y_{v_2}<\tilde{Y}_{i,j}
\le y_{v_2-1}+t_{v_3-1}} \quad\mbox{and} \quad\overline{\xi}_{i,j}(v_2,v_3)
= \mathbf{1}_{y_{v_2-1}<\tilde
{Y}_{i,j} \le y_{v_2}+t_{v_3}}.
\end{eqnarray*}
Clearly, for any $(x,y,\delta)\in I_{v_1,v_2,v_3}$, we have
$\underline{\xi}_{i,j}(v_2,v_3) \le\tilde\xi_{i,j}(\delta,y)
\le \overline{\xi}_{i,j}(v_2,v_3)$. Since $N_n\to\infty$ and $\delta
_n\to0$,
there exists a constant $c_2<\infty$ such that
$0\le\E[\overline{\xi}_{i,j}(v_2,v_3)] - \E[
\underline{\xi}_{i,j}(v_2,v_3) ] \le c_2 N_n^{1/q} /\ell_n $.
Additionally, for $x\in[x_{v_1-1},x_{v_1}]$, by Condition \ref{assump:r2},\
$|\varpi_{i,j}(x)-\varpi_{i,j}(x_{v_1})|\le\tau_n |x-x_{v_1}|\le\tau
_n b/\ell_n$. Thus, there exists a constant $c_3<\infty$ such that
%
%
\begin{eqnarray}
\label{eq:upper} && \varpi_{i,j}(x) \bigl\{ \tilde\xi_{i,j}(
\delta,y) - \E\bigl[ \tilde\xi_{i,j}(\delta,y) \bigr]\bigr\}
\nonumber\\
&&\quad \le
\varpi_{i,j}(x_{v_1}) \bigl\{ \overline{\xi}_{i,j}(v_2,v_3)
- \E\bigl[ \underline{\xi}_{i,j}(v_2,v_3) \bigr]
\bigr\} + \tau_n b/\ell_n
\\
& &\quad\le
\varpi_{i,j}(x_{v_1}) \bigl\{ \overline{\xi}_{i,j}(v_2,v_3)
- \E\bigl[ \overline{\xi}_{i,j}(v_2,v_3) \bigr]
\bigr\} + c_3 \bigl(\tau_n +N_n^{1/q}
\bigr)/\ell _n,\nonumber
\end{eqnarray}
uniformly over $i,j$, and $(x,y,\delta)\in I_{v_1,v_2,v_3}$. Similarly,
%
%
\begin{eqnarray}
\label{eq:lower} && \varpi_{i,j}(x) \bigl\{ \tilde\xi_{i,j}(
\delta,y) - \E\bigl[ \tilde\xi_{i,j}(\delta,y) \bigr]\bigr\}
\nonumber
\\[-8pt]
\\[-8pt]
\nonumber
&&\quad \ge
\varpi_{i,j}(x_{v_1}) \bigl\{ \underline{\xi}_{i,j}(v_2,v_3)
- \E\bigl[ \underline{\xi}_{i,j}(v_2,v_3) \bigr]
\bigr\} - c_3 \bigl(\tau_n +N_n^{1/q}
\bigr)/\ell_n.
\end{eqnarray}
Combining (\ref{eq:upper}) and (\ref{eq:lower}) and using $N_n(\tau_n
+N_n^{1/q})/\ell_n=\mathrm{O}(1)$, we have
%
%
\begin{equation}
\label{eq:upperlower} \sup_{x,y,\delta} \bigl|\tilde{D}_n(
\delta,x,y)\bigr| \le\max_{v\in\mathcal
{V}_n} \bigl\{ \bigl|\underline{
\Delta}_n(v)\bigr| + \bigl|\overline{\Delta}_n(v)\bigr| \bigr\} +\mathrm{O}(1),
\end{equation}
where $v=(v_1,v_2,v_3)$,
\begin{eqnarray*}
\underline{\Delta}_n(v) &=& \sum^n_{i=1}
\sum^{m_i}_{j=1} \varpi_{i,j}(x_{v_1})
\bigl\{ \underline {\xi}_{i,j}(v_2,v_3) - \E
\bigl[ \underline{\xi}_{i,j}(v_2,v_3) \bigr]
\bigr\},
\\
\overline{\Delta}_n(v) &=& \sum
^n_{i=1} \sum^{m_i}_{j=1}
\varpi_{i,j}(x_{v_1}) \bigl\{ \overline{\xi
}_{i,j}(v_2,v_3) - \E\bigl[ \overline{
\xi}_{i,j}(v_2,v_3) \bigr] \bigr\}.
\end{eqnarray*}
We now apply Theorem \ref{thmm:cp} to $\underline{\Delta}_n(v)$ and
$\overline{\Delta}_n(v)$.
For $\chi_n$ in (\ref{eq:thmcpcon}), with $\phi_n$ in (\ref{eq:conr2})
and $\E[\overline{\xi}_{i,j}(h_2,h_3)]=\mathrm{O}(\delta_n+N_n^{1/q}/\ell
_n)=\mathrm{O}(\delta_n)$,
we can take $\chi_n=\mathrm{O}(\delta_n\phi_n)$. By Theorem \ref{thmm:cp}(ii),
$\max_{v\in\mathcal{V}_n} |\overline{\Delta}_n(v)|=\mathrm{O}_\mathrm{p}\{[\delta
_n\phi_n (\log N_n)^3]^{1/2}\}$.
The latter bound also holds for $ \max_{v\in\mathcal{V}_n} |\underline
{\Delta}_n(v)|$. The desired result then follows from (\ref{eq:upperlower}).
\end{pf*}

\subsection{Asymptotic expansions}\label{sec:expansion}
Throughout the proofs, we use the following notation:
\begin{eqnarray*}
L_{\mu}(\delta_1,x) &=& \sum
^n_{i=1} \sum^{m_i}_{j=1}
K_{b_n}(x_{i,j}-x) \mathbf{1}_{Y_{i,j}\le\mu(x)+\delta_1},
\\
L_{\mu}(x) &=& \sum^{n}_{i=1}
\sum^{m_i}_{j=1} K_{b_n}(x_{i,j}-x),
\\
J_{\mu}(\delta_1,x) &=& \E\bigl[ L_{\mu}(
\delta_1,x) \bigr],
\\
L_{s} (\delta_1,\delta_2,x) &=& \sum
^{n}_{i=1} \sum
^{m_i}_{j=1} K_{h_n}(x_{i,j}-x)
\mathbf{1}_{ |Y_{i,j}-\mu(x)-\delta_1|\le s(x) + \delta_2 },
\\
L_{s}(x) &=& \sum^{n}_{i=1}
\sum^{m_i}_{j=1} K_{h_n}
(x_{i,j}-x),
\\
J_{s}(\delta_1,\delta_2,x)&=&\E\bigl[
L_{s}(\delta_1,\delta_2,x) \bigr].
\end{eqnarray*}

\begin{lemma}\label{lem:1}
Assume that Conditions
\ref{con:reg}--\ref{con:lc} hold. Then, we have
\begin{longlist}[(ii)]
\item[(i)] Uniformly over $x\in{\cal S}_\epsilon[a,b]$,
%
%
\begin{equation}
\label{eq:lem1a} \sum^n_{i=1} \sum
^{m_i}_{j=1} \biggl( \frac{x_{i,j}-x}{b_n}
\biggr)^r K \biggl( \frac{x_{i,j}-x}{b_n} \biggr) = \frac{N_n b_n}{b-a}
\int_{\R} u^r K(u)\,\mathrm{d}u + \mathrm{O}(1).
\end{equation}
\item[(ii)] Let $g(x,v)$ be a measurable bivariate function on
$[a,b]^2$. Define
%
%
\begin{equation}
{\cal G}_g(x) = \sum^{n}_{i=1}
\sum^{n_i}_{j=1} g(x,x_{i,j})
K_{b_n}(x_{i,j}-x).
\end{equation}
Further assume that $\sup_{x\in[a,b]}|\partial^{s}
(x,v)/\partial v^s|<\infty, s=0,1,\ldots,r$ for some $r\in\N$.
Then uniformly over $x\in{\cal S}_\epsilon[a,b]$,
%
%
\begin{equation}
\label{eq:lem1b} {\cal G}_g(x)= \sum^{r-1}_{s=0}
\frac{\partial^s g(x,v)}{\partial v^s} \Big|_{v=x} \frac{N_n b^{s+1}_n}{(b-a)s!} \int_{\R}
u^s K(u)\,\mathrm{d}u + \mathrm{O}\bigl(1+N_n b^{r+1}_n
\bigr).
\end{equation}
\end{longlist}
\end{lemma}

\begin{pf}
(i) Recall the ordered locations $\tilde{x}_k$ in Condition
\ref{con:reg}. Define
%
%
\begin{eqnarray}
\label{eq:snx} S_n(x) &=& \sum^{N_n}_{k=1}
\biggl( \frac{\tilde{x}_k-x}{b_n} \biggr)^r K \biggl( \frac{\tilde{x}_k-x}{b_n}
\biggr),\label{eq:snx}
\\
I_n(x) &=& \sum^{N_n}_{k=0} (
\tilde{x}_{k+1}-\tilde{x}_k) \biggl( \frac{\tilde{x}_k-x}{b_n}
\biggr)^r K \biggl( \frac{\tilde
{x}_k-x}{b_n} \biggr), \label{eq:inx}
\\
\varrho_n &=& \max_{0\le k\le N_n} \bigl|\tilde{x}_{k+1}-
\tilde{x}_k - (b-a)/N_n \bigr| =\mathrm{O}\bigl(N^{-2}_n
\bigr), \label{eq:ellnx}
\\
{\cal I}(x) &=& \bigl\{1\le k\le N_n\dvt \tilde{x}_k-x
\in\bigl[- b_n - (b-a)/N_n - \varrho_n,
b_n \bigr] \bigr\}. \label{eq:ix}
\end{eqnarray}

Assume without loss of generality that $K$ has
support $[-1,1]$. Condition (\ref{eq:dense}) implies that $\sup_{x\in
[a,b]}|{\cal I}(x)|=\mathrm{O}(N_nb_n)$, where and hereafter $|{\cal
I}|$ is the cardinality of a set ${\cal I}$. Because $K$ has
support $[-1,1]$, $K_{b_n}(\tilde{x}_k-x)=0$ for $k\notin
{\cal I}(x)$. Additionally, for $k\in{\cal I}(x)$, the summands
in $S_n(x)$ are uniformly bounded. Thus,
%
%
\begin{eqnarray}
\label{eq:i2bound} S_n(x) = \sum_{k\in{\cal I}(x)}
\biggl( \frac{\tilde{x}_k-x}{b_n} \biggr)^r K \biggl( \frac{\tilde
{x}_k-x}{b_n}
\biggr) = \mathrm{O}\bigl[\bigl|{\cal I}(x)\bigr|\bigr] = \mathrm{O}(N_nb_n),
\end{eqnarray}
uniformly over $x\in[a,b]$.

By (\ref{eq:dense}), elementary calculation shows that,
uniformly over $x\in{\cal S}_\epsilon[a,b]$,
%
%
\begin{eqnarray}
\label{eq:snapp} \frac{b-a}{N_n}S_n(x) - I_n(x) &=& -
\sum^{N_n}_{k=1} \biggl(\tilde{x}_{k+1}-
\tilde{x}_k - \frac{b-a}{N_n} \biggr) \biggl( \frac{\tilde{x}_k-x}{b_n}
\biggr)^r K \biggl( \frac{\tilde
{x}_k-x}{b_n} \biggr)
\nonumber
\\[-8pt]
\\[-8pt]
\nonumber
&=& \mathrm{O}(
\varrho_n) \sup_{x\in[a,b]} \bigl|S_n(x)\bigr| =
\mathrm{O}(b_n/N_n).
\end{eqnarray}

Write $u_k=(\tilde{x}_k-x)/b_n$.
Observe that $I_n(x)=\sum^{N_n}_{k=0} \int^{\tilde{x}_{k+1}}_{\tilde{x}_k}
u_k^r K (u_k)\,\mathrm{d}v$. Thus, by the
triangle inequality, we have
%
%
\begin{eqnarray}\label{eq:VKa}
\biggl\llvert I_n(x) - \int^{\tilde{x}_{N_n+1}}_{\tilde{x}_0}
\biggl(\frac{v-x}{b_n} \biggr)^r K \biggl(\frac{v-x}{b_n}
\biggr)\,\mathrm{d}v \biggr\rrvert \le \sum^{N_n}_{k=0}
V_k,
\nonumber
\\[-8pt]
\\[-8pt]
\eqntext{\mbox{where } V_k = \displaystyle\int^{\tilde{x}_{k+1}}_{\tilde{x}_k}
\biggl\llvert u^r_k K(u_k)- \biggl(
\frac{v-x}{b_n} \biggr)^r K \biggl(\frac{v-x}{b_n} \biggr)
\biggr\rrvert \,\mathrm{d}v.\qquad\qquad}
\end{eqnarray}
Since $K$ has bounded derivative, $|y^r K(y)-z^r K(z)|=\mathrm{O}(|y-z|)$ for
$y,z\in[-1,1]$.
Also, $|u_k-(v-x)/b_n| =|v-\tilde{x}_k|/b_n$. Thus,
under Condition \ref{con:reg},
%
%
\begin{equation}
\label{eq:VKb} |V_k|= \mathrm{O}(1)\int^{\tilde{x}_{k+1}}_{\tilde{x}_k}
\frac{v-\tilde{x}_k}{b_n}\,\mathrm{d}v = \frac{\mathrm{O}[(\tilde{x}_{k+1}-\tilde{x}_{k})^2]}{b_n} = \frac{\mathrm{O}(1)}{N^2_n b_n}.
\end{equation}
Furthermore, it is easily seen that, for $k\notin{\cal I}(x)$,
$\min(|\tilde{x}_k-x|,|\tilde{x}_{k+1}-x|)> b_n$, which implies
$K(u_k)=0, K\{(v-x)/b_n\}=0$ for $v\in[\tilde{x}_k,\tilde{x}_{k+1}]$,
and consequently $V_k=0$.
Thus, by (\ref{eq:VKa}) and~(\ref{eq:VKb}),
%
%
\begin{equation}
\label{eq:experr} \biggl\llvert I_n(x) - \int^{\tilde{x}_{N_n+1}}_{\tilde{x}_0}
\biggl(\frac{v-x}{b_n} \biggr)^r K \biggl(\frac{v-x}{b_n}
\biggr)\,\mathrm{d}v \biggr\rrvert \le \sum_{k\in{\cal I}(x)}
V_k = \mathrm{O}(1/N_n),
\end{equation}
uniformly over $x\in{\cal S}_\epsilon[a,b]$,

Notice that $\sum^n_{i=1} \sum^{m_i}_{j=1} [(x_{i,j}-x)/b_n]^r
K_{b_n}(x_{i,j}-x )=S_n(x)$. Recall that $\tilde{x}_0=a$ and
$\tilde{x}_{N_n+1}=b$. The desired result then follows from
(\ref{eq:snapp}) and (\ref{eq:experr}) in view of
\begin{eqnarray*}
\int^{\tilde{x}_{N_n+1}}_{\tilde{x}_0} \biggl(\frac{v-x}{b_n}
\biggr)^r K \biggl(\frac{v-x}{b_n} \biggr)\,\mathrm{d}v =b_n \int
^{(b-x)/b_n}_{(a-x)/b_n} u^rK(u)\,\mathrm{d}u = b_n
\int^1_{-1}u^rK(u)\,\mathrm{d}u
\end{eqnarray*}
for all $x\in{\cal S}_\epsilon[a,b]$ and large enough $n$.

(ii) The expression (\ref{eq:lem1b}) easily follows from (i)
in view of the Taylor expansion $g(x,x_{i,j})=\sum^{r-1}_{s=0}
\partial^s g(x,v)/\partial v^s|_{v=x}
(x_{i,j}-x)^s/s!+\mathrm{O}(b^r_n)$ for $|x_{i,j}-x|\le b_n$.
\end{pf}

\begin{lemma}\label{lem:expansion}
Assume that Conditions \ref{con:reg}--\ref{con:lc} hold. Let
$\rho_\mu(x),\rho_s(x),\kappa,\kappa_+$ be as in Theorems~\ref{thmm:m}--\ref{thmm:v}. Then, for $\delta_1\to
0,\delta_2\to0$, we have uniformly over $x\in{\cal
S}_\epsilon[a,b]$,
\begin{eqnarray*}
J_\mu(0,x) &=& L_\mu(x)/2 -N_nb^3_n
\rho_\mu(x) f_e(0) \psi _K/\bigl[(b-a)s(x)
\bigr]+\mathrm{O}\bigl(1+N_nb^5_n\bigr),
\label{eq:jsx2}
\\
J_\mu(\delta_1,x) &=& J_\mu(0,x) +
N_n b_n \delta_1 \bigl\{
f_e(0)/\bigl[(b-a)s(x)\bigr] + \mathrm{O}\bigl[(N_nb_n)^{-1}+b^2_n+
\delta_1\bigr] \bigr\}, \label {eq:jsxdiff}
\\
J_s(\delta_1,0,x) &=& L_s(x)/2 -
N_n h_n\kappa_{+} \bigl\{ \bigl[
{h_n^2} \psi_K \rho_s(x) -
\delta_1 \kappa\bigr] /\bigl[(b-a)s(x)\bigr] + \mathrm{O}\bigl(h_n^4
+ \delta_1^2\bigr) \bigr\},
\\
J_s(\delta_1,\delta_2,x) &=&
J_s(\delta_1,0,x) + N_n h_n
\delta_2 \bigl\{ \kappa_{+} /\bigl[(b-a) s(x)\bigr] + \mathrm{O}
\bigl( h_n^2 + \delta_1 +
\delta_2\bigr) \bigr\}.
\end{eqnarray*}
\end{lemma}

\begin{pf}
Recall that $F_e$ and $f_e$ are the distribution and density
functions of $e_{i,j}$. The assumption ${\cal Q}(e_{i,j})=0$
implies that $F_e(0)=1/2$. Notice that
\begin{eqnarray*}
J_\mu(0,x)-L_\mu(x)/2 &=& \sum
^n_{i=1} \sum^{m_i}_{j=1}
K_{b_n}(x_{i,j}-x) \bigl[ \p\bigl\{Y_{i,j}\le\mu(x)
\bigr\} - 1/2 \bigr]
\\
&=& \sum^n_{i=1}
\sum^{m_i}_{j=1} K_{b_n}(x_{i,j}-x)
g(x,x_{i,j}),
\end{eqnarray*}
where $g(x,v)= F_e \{ [\mu(x)-\mu(v)]/s(v) \} -
F_e(0)$. The symmetry of $K$ entails
$\int u^s K(u)\,\mathrm{d}u=0, s=1,3$. The first expression
then follows from Lemma \ref{lem:1}(ii) with $r=4$.

Similarly, we can show $J'_\mu(0,x):=\partial
J_\mu(\delta_1,x)/\partial\delta_1|_{\delta_1=0}=N_n b_n
f_e(0)/[(b-a)s(x)]+\mathrm{O}(1+N_nb^3_n)$ and $J''_\mu(\delta_1,x):=\partial^2
J_\mu(\delta_1,x)/\partial\delta_1^2=\mathrm{O}(Nb_n)$ uniformly over $\delta_1,x$.
So, the second expression
follows from the Taylor expansion
$J_\mu(\delta_1,x)-J_\mu(0,x)=\delta_1 J'_\mu(0,x)+\mathrm{O}(Nb_n\delta_1^2)$.
The other two expressions can be similarly treated.
We omit the details.
\end{pf}

\subsection{Proof of Theorems \texorpdfstring{\protect\ref{thmm:m}--\protect\ref{thmm:v}}{4.1--4.2}}\label{sec:pfthm}

Let
$L_\mu(x),L_\mu(\delta_1,x),J_\mu(\delta_1,x),L_s(x),L_s(\delta_1,\delta_2,x)$
and $J_s(\delta_1,\delta_2,x)$ be as in Section~\ref{sec:expansion}.

\begin{pf*}{Proof of Theorem \ref{thmm:m}}
Let $\delta_n=[(\log N_n)^3/(N_nb_n)]^{1/2} + b_n^2\to0$. Let
$l_n \uparrow\infty$ be a positive sequence satisfying
$\delta_n l_n \to0$. First, we show $\hat\Delta_\mu(x) :=
\hat{\mu}(x) - \mu(x) = \mathrm{O}_\mathrm{p}(l_n \delta_n)$ uniformly over
$x\in{\cal S}_\epsilon([a,b])$. Since $\hat{\mu}(x)$ is a
solution to (\ref{eq:muest}), by
Koenker (\cite{koe2005},\break  pages~32--33),
%
%
\begin{equation}
\label{eq:koenker} \bigl|L_\mu\bigl(\hat\Delta_\mu(x),x\bigr)
- L_\mu(x)/2\bigr| \le \sum_{i,j}
K_{b_n}(x_{i,j}-x) \mathbf{1}_{Y_{i,j}=\hat\mu(x)} =
\mathrm{O}_\mathrm{p}(1),
\end{equation}
uniformly over $x$.
Let
\[
\Omega_n(x) = \bigl[L_\mu(l_n
\delta_n,x) - J_\mu(l_n \delta_n,x)
\bigr] - \bigl[L_\mu(0,x)-J_\mu(0,x)\bigr].
\]
We can apply Theorem \ref{pro:osc} with $\varpi
_{i,j}(x)=K_{b_n}(x_{i,j}-x)$ to $\Omega_n(x)$. For $\tau_n$ and $\phi
_n$ in\break Condition~\ref{assump:r2}, $\tau_n=\mathrm{O}(1/b_n)$ and $\phi_n=\mathrm{O}(N_n b_n)$ (see Lemma
\ref{lem:1}). By Theorem \ref{pro:osc},\break $\sup_{x\in[a,b]} |\Omega_n(x)|=
\mathrm{O}_\mathrm{p}\{[N_nb_n l_n\delta_n(\log N_n)^3]^{1/2}\}$.
By the same argument, we can show
%
%
\begin{equation}
\label{eq:osc2} \sup_{x\in[a,b]} \bigl|L_\mu(0,x)-J_\mu(0,x)\bigr|
= \mathrm{O}_\mathrm{p}\bigl\{\bigl[N_nb_n (\log
N_n)^3\bigr]^{1/2}\bigr\}.
\end{equation}
Hence, by (\ref{eq:osc2}) and Lemma
\ref{lem:expansion}, uniformly over $x\in{\cal
S}_\epsilon([a,b])$,
%
%
\begin{eqnarray}
\label{eq:consis} L_\mu(l_n \delta_n,x) -
L_\mu(x)/2 &=&\bigl[ J_\mu(l_n
\delta_n,x) - J_\mu(0,x) \bigr] + \bigl[
J_\mu(0,x) - L_\mu(x)/2 \bigr]
\nonumber\\
&&{} + \bigl[
L_\mu(0,x) - J_\mu(0,x) \bigr] + \Omega_n(x)
\\
& = & N_nb_n l_n \delta_n
f_e(0)/\bigl[(b-a)s(x)\bigr] \bigl[ 1 + \mathrm{o}(1) \bigr] +
\mathrm{O}_\mathrm{p}(\nu_n),\nonumber
\end{eqnarray}
where $\nu_n= N_n b^3_n + 1 + [N_nb_n(\log N_n)^3]^{1/2} +
[N_nb_nl_n\delta_n(\log N_n)^3]^{1/2}$. Because $l_n\to\infty$
and $l_n\delta_n\to0$, it is easy to see that $\nu_n=\mathrm{o}(N_nb_n
l_n \delta_n)$ and $N_nb_n l_n \delta_n\to\infty$, which
implies $L_\mu(l_n \delta_n,x) - L_\mu(x)/2\to\infty$ uniformly
over $x\in{\cal S}_\epsilon[a,b]$ in view of
$\sup_{x}s(x)<\infty$. Since $L_\mu(\delta_1,x)$ is
nondecreasing in $\delta_1$, (\ref{eq:koenker}) and
(\ref{eq:consis}) entail $\p\{\sup_x \hat\Delta_\mu(x) \le l_n
\delta_n \}\to1$. Similarly, $\p\{\inf_x \hat\Delta_\mu(x)
\ge-l_n \delta_n \}\to1$. So, $\sup_x
|\hat\Delta_\mu(x)|=\mathrm{O}_\mathrm{p}(l_n \delta_n)$. Since the rate of
$l_n \to\infty$ can be arbitrarily slow, $\sup_x
|\hat\Delta_\mu(x)|=\mathrm{O}_\mathrm{p}( \delta_n)$.

Again, by (\ref{eq:koenker}) and
Lemma \ref{lem:expansion}, uniformly over $x\in{\cal
S}_\epsilon([a,b])$,
\begin{eqnarray*}
L_\mu(0,x)-J_\mu(0,x)&=& L_\mu\bigl(\hat
\Delta_\mu(x),x\bigr)-J_\mu\bigl(\hat\Delta_\mu
(x),x\bigr) + \mathrm{O}_\mathrm{p} \bigl[ \sqrt{N_nb_n
\delta_n (\log N_n)^3 } \bigr]
\\
&=&
\bigl[L_\mu\bigl(\hat\Delta_\mu(x),x\bigr) -
L_\mu(x)/2\bigr] + \bigl[ L_\mu(x)/2-J_\mu(0,x)
\bigr]
\\
&&{} - \bigl[J_\mu\bigl(\hat\Delta_\mu(x),x
\bigr) - J_\mu(0,x) \bigr] + \mathrm{O}_\mathrm{p} \bigl[
\sqrt{N_nb_n \delta_n (\log
N_n)^3 } \bigr]
\\
&=& \mathrm{O}_\mathrm{p}(1) +
N_n b^3_n \rho_\mu(x)f_e(0)
\psi_K/\bigl[(b-a)s(x)\bigr] + \mathrm{O}\bigl(1+N_nb^5_n
\bigr)
\\
&&{} - N_nb_n \hat\Delta_\mu(x) \bigl
\{ f_e(0)/\bigl[(b-a)s(x)\bigr] + \mathrm{O}(\delta_n)\bigr\}
\\
&&{} +\mathrm{O}_\mathrm{p} \bigl[ \sqrt{N_nb_n
\delta_n (\log N_n)^3 } \bigr].
\end{eqnarray*}
The representation (\ref{eq:mbahadur}) then follows by solving
$\hat\Delta_\mu(x)$ from the above equation.
\end{pf*}

\begin{pf*}{Proof of Theorem \ref{thmm:v}}
We use the argument in Theorem \ref{thmm:m} and only sketch the outline.
Let
\[
D_s(\delta_1,\delta_2,x)=
\bigl[L_s(\delta_1,\delta_2,x)-J_s(
\delta _1,\delta_2,x)\bigr]- \bigl[L_s(0,0,x)-J_s(0,0,x)
\bigr].
\]
Using Theorem \ref{pro:osc}, we can show that
%
%
\begin{eqnarray}
\sup_{x\in[a,b]}\bigl|L_s(0,0,x)-J_s(0,0,x)\bigr|&=&\mathrm{O}_\mathrm{p}
\bigl\{\bigl[N_nh_n (\log N_n)^3
\bigr]^{1/2}\bigr\}, \label{eq:oscaa}
\\
\sup_{|\delta_1|+|\delta_2|\le\delta_n, x\in[a,b]} \bigl|D_s(\delta _1,
\delta_2,x)\bigr|&=&\mathrm{O}_\mathrm{p}\bigl\{\bigl[N_nh_n
(\log N_n)^3\bigr]^{1/2}\bigr\}, \label{eq:oscs}
\end{eqnarray}
hold for all $b_n \to0, h_n\to0$ and $\delta_n\to0$
satisfying $\sup_{n}\log N_n/[N_n \min(b_n,h_n) \delta_n]<\infty$.

Let $\delta_n=b^2_n + h^2_n + [(\log N_n)^3/(N_nb_n)]^{1/2}
+[(\log N_n)^3/(N_nh_n)]^{1/2}$ and $l_n\to\infty$ be a
sequence such that $l_n\delta_n\to0$. By Theorem \ref{thmm:m},
$ \tilde\Delta_\mu(x):=\tilde{\mu}(x) - \mu(x) = \mathrm{O}_\mathrm{p}
(\delta_n )$. Using (\ref{eq:oscs}) and Lemma~\ref{lem:expansion}, we can derive the following counterpart of
(\ref{eq:consis})
\begin{eqnarray*}
L_s\bigl(\tilde\Delta_\mu(x), l_n
\delta_n, x\bigr) - L_s(x)/2 &=& \bigl[J_s
\bigl(\tilde\Delta_\mu(x), l_n \delta_n,x
\bigr) -J_s\bigl(\tilde\Delta_\mu (x),0,x\bigr)\bigr]
\\
& & {}+ \bigl[J_s\bigl(\tilde\Delta_\mu(x),0,x
\bigr)-L_s(x)/2\bigr] + L_s(0,0,x)-J_s(0,0,x)
\\
&& {} + \mathrm{O}_\mathrm{p} \bigl\{\bigl[ N_nh_n
l_n \delta_n (\log N_n )^3
\bigr]^{1/2} \bigr\}
\\
&=& N_nh_n l_n \delta_n
\kappa_{+}/\bigl[(b-a)s(x)\bigr] \bigl[ 1 + \mathrm{o}_\mathrm{p}(1)
\bigr] \to\infty.
\end{eqnarray*}
Let $\hat\Delta_s(x)=\hat{s}(x)-s(x)$. By the same argument in
(\ref{eq:koenker}), $\sup_{x}|L_s(\tilde\Delta_\mu(x),
\hat\Delta_s(x), x)-L_s(x)/2|=\mathrm{O}_\mathrm{p}(1)$. Notice that
$L_s(\tilde\Delta_\mu(x),\delta_2,x)$ is nondecreasing in $x$.
Thus, $\p\{\sup_x \hat\Delta_s(x) \le l_n k_n\}\to1$.
Similarly, $\p\{\inf_x \hat\Delta_s(x) \ge-l_n k_n\}\to1$.
Then $\sup_x |\hat\Delta_s(x)|=\mathrm{O}_\mathrm{p}(\delta_n)$.

Write $\varpi_n=[ N_nh_n \delta_n (\log N_n )^3]^{1/2}$. To
derive the Bahadur representation (\ref{eq:vbahadur}), we use
(\ref{eq:oscs}) and Lemma \ref{lem:expansion} to obtain
\begin{eqnarray*}
&&L_s(0,0,x)-J_s(0,0,x) \\
&&\quad= \bigl[L_s\bigl(
\tilde\Delta_\mu(x),\hat\Delta_s(x),x
\bigr)-L_s(x)/2\bigr] + \bigl[L_s(x)/2-J_s
\bigl(\hat\Delta_\mu(x),0,x\bigr)\bigr]
\\
&&\qquad{} -\bigl[J_s
\bigl(\tilde\Delta_\mu(x),\hat\Delta_s(x),x
\bigr)-J_s\bigl(\tilde\Delta_\mu(x),0,x\bigr)\bigr] +
\mathrm{\mathrm{O}}_\mathrm{p} (\varpi_n)
\\
&&\quad= \mathrm{O}_\mathrm{p}(1) +
N_n h_n \kappa_{+} \bigl\{
\bigl[h^2_n\psi_K \rho_s(x) -
\kappa\tilde\Delta_u(x)\bigr]/\bigl[(b-a)s(x)\bigr] +\mathrm{O}
\bigl(h^4_n + \delta^2_n\bigr)
\bigr\}
\\
&&\qquad{} - N_n h_n \hat\Delta_s(x)
\bigl\{\kappa_{+}/\bigl[(b-a)s(x)\bigr]+\mathrm{O}(\delta_n)\bigr\}
+ \mathrm{O}_\mathrm{p} (\varpi_n).
\end{eqnarray*}
Solving $\hat\Delta_s(x)$ from the above equation, we obtain
the Bahadur representation (\ref{eq:vbahadur}).
\end{pf*}

\subsection{Proof of Corollaries \texorpdfstring{\protect\ref{cor:m}--\protect\ref{cor:v}}{4.1--4.2}} \label{sec:clt}
Again we use the coupling argument to convert
the dependent data to $m$-dependent case. Theorem \ref{thmm:romanowolf}
below presented a CLT for $m$-dependent sequence with unbounded $m$.

\begin{thmm}[(Romano and Wolf \cite{rom2000})]\label{thmm:romanowolf}
Let $Z_{n,j},1\le j\le d_n$, be a triangular array of mean zero
$k_n$-dependent random variables. Define
\begin{eqnarray*}
S_n=\sum^{d_n}_{j=1}
Z_{n,j},\qquad B^2_n = \operatorname{Var}
(S_n),\qquad S_{n,h,a}=\sum^{a+h-1}_{j=a}
Z_{n,j}, \qquad B^2_{n,h,a}=\operatorname{Var}
(S_{n,h,a}).
\end{eqnarray*}
Assume that there exist some $\delta>0, -1\le\gamma<1,
C_{n,1},C_{n,2},C_{n,3}>0$ such that
\begin{eqnarray*}
\begin{array} {l@{\qquad}l} \textup{(a)}\quad \E\bigl(|Z_{n,j}|^{2+\delta}
\bigr) =\mathrm{O}(C_{n,1}); &\textup{(b)}\quad B^2_{n,h,a}/h^{1+\gamma}
=\mathrm{O}(C_{n,2})\qquad \mbox{for all } h\ge k_n,a;
\\
\textup{(c)}\quad B^2_n/\bigl(d_n C_{n,2}^\gamma
\bigr)\ge C_{n,3}; &\mathrm{(d)}\quad C_{n,2}/C_{n,3}=\mathrm{O}(1);
\\
\textup{(e)}\quad C_{n,1}/C_{n,3}^{(2+\delta)/2}=\mathrm{O}(1); &\mathrm{(f)}\quad
k_n^{1+(1-\gamma)(1+2/\delta)}/d_n\to0. \end{array} %
\end{eqnarray*}
Then $S_n/B_n\Rightarrow N(0,1)$.
\end{thmm}

\begin{pf*}{Proof of Corollaries \ref{cor:m}--\ref{cor:v}}
We only prove Corollary \ref{cor:m} since Corollary \ref{cor:v} can
be similarly treated. By the Bahadur representation
(\ref{eq:mbahadur}), under the specified condition,
$r_n\sqrt{N_nb_n}\to0$. Thus, it suffices to show
$(N_nb_n)^{-1/2} Q_{b_n}(x) \Rightarrow
N(0,\varphi_K/[4(b-a)])$. Recall $e_{i,j}(k_n)$ and
$\tilde{Y}_{i,j}$ in (\ref{eq:couple}) and
(\ref{eq:Fn2}). Define the coupling process
\[
\tilde{Q}_{b_n}(x) = - \sum^{n}_{i=1}
\sum^{m_i}_{j=1} \bigl\{
\mathbf{1}_{\tilde{Y}_{i,j}\le\mu(x)} - \E[\mathbf{1}_{\tilde
{Y}_{i,j}\le\mu(x)}]\bigr\}
K_{b_n}(x_{i,j}-x).
\]
Let the coupling lag $k_n=\lfloor c\log N_n \rfloor$ be chosen
as in Theorem \ref{lem:appr}. By Theorem \ref{lem:appr},
$Q_{b_n}(x)-\tilde{Q}_{b_n}(x)=\mathrm{O}_\mathrm{p}[(\log N_n)^2]=\mathrm{o}_\mathrm
{p}[(N_nb_n)^{1/2}]$. It remains to show $(N_nb_n)^{-1/2}
\tilde{Q}_{b_n}(x) \Rightarrow \break N(0,\varphi_K/[4(b-a)])$. Recall
$M_n=\max_{1\le i\le n} m_i$. Set $\tilde{Y}_{i,j}=0$ for
$m_i<j\le M_n$. Define
\begin{eqnarray*}
Z_{n,j}=\sum^n_{i=1}
\zeta_{i,j}, \qquad\mbox{where } \zeta_{i,j} = (N_nb_n)^{-1/2}
\bigl\{ \mathbf{1}_{\tilde{Y}_{i,j}\le\mu(x)} - \E[\mathbf{1}_{\tilde
{Y}_{i,j}\le\mu(x)}]\bigr\}
K_{b_n}(x_{i,j}-x).
\end{eqnarray*}
Then we can write $-(N_nb_n)^{-1/2} \tilde{Q}_{b_n}(x) =
\sum^{M_n}_{j=1} Z_{n,j}$. Notice that $Z_{n,j},j=1,2,\ldots,$
are $(2k_n+1)$-dependent, and
$\zeta_{i,j},i=1,2,\ldots,$ are independent
for each fixed $j$.

Let $S_n,B^2_n,S_{n,h,a}$ and $B^2_{n,h,a}$ be defined in
Theorem \ref{thmm:romanowolf}. We shall verify the conditions in
Theorem \ref{thmm:romanowolf}. By
the independence of the summands $\zeta_{i,j}$ in $Z_{n,j}$,
\begin{eqnarray*}
\E\bigl(|Z_{n,j}|^4\bigr) &=& \sum
^n_{i=1} \E\bigl(|\zeta_{i,j}|^4
\bigr) + 6 \sum_{i_1\ne i_2} \E\bigl(|\zeta_{i_1,j}|^2
\bigr) \E\bigl(|\zeta_{i_2,j}|^2\bigr)
\\
&=&
\frac{\mathrm{O}(1)}{(N_nb_n)^2} \Biggl\{ \sum^n_{i=1}
K^4_{b_n}(x_{i,j}-x) + \Biggl[ \sum
^n_{i=1} K^2_{b_n}(x_{i,j}-x)
\Biggr]^2 \Biggr\} = \mathrm{O} \bigl( 1/M_n^2\bigr),
\end{eqnarray*}
in view of $nM_n=\mathrm{O}(N_n)$. Since $\tilde{Y}_{i,j}$ and $Y_{i,j}$
have same distribution, we have $g(x,x_{i,j}):=\operatorname{Var}({\bf
1}_{\tilde{Y}_{i,j}\le\mu(x)}) =
F_e\{[\mu(x)-\mu(x_{i,j})]/s(x_{i,j})\}-F^2_e\{[\mu(x)-\mu
(x_{i,j})]/s(x_{i,j})\}$.
Recall $F_e(0)=1/2$. Then $g(x,x)=1/4$. Thus, by
(\ref{eq:conclt}) and the $(2k_n+1)$-dependence of
$\tilde{Y}_{i,j},j\in\Z$, applying Lemma \ref{lem:expansion}(ii) with $r=1$ produces
\begin{eqnarray*}
B^2_n &=& \frac{1}{N_nb_n} \sum
^n_{i=1} \sum^{n_i}_{j=1}
\operatorname{Var}(\mathbf{1}_{\tilde{Y}_{i,j}\le\mu(x)}) K^2_{b_n}(x_{i,j}-x)
\\
&&{} + \frac{\mathrm{O}(1)}{N_nb_n} \sum^n_{i=1}
\sum_{1\le j_1<j_2 \le n_i,
|j_1-j_2|\le2k_n} K_{b_n}(x_{i,j_1}-x)K_{b_n}(x_{i,j_2}-x)
\\
&=& \frac{1}{N_nb_n} \biggl[\frac{N_nb_n \varphi_K}{4(b-a)} + \mathrm{O}\bigl(N_nb^2_n
\bigr) \biggr] + \frac{\mathrm{O}(nM_n k_n b_n \iota_n)}{N_nb_n} \to\frac
{\varphi_K}{4(b-a)},
\end{eqnarray*}
in view of $nM_n=\mathrm{O}(N_n)$ and $k_n\iota_n\to0$.
Similarly, we can show $B^2_{n,h,a}=\mathrm{O}(nh/N_n)=\mathrm{O}(h/M_n)$.
Therefore, it is easy to see that the conditions in Theorem
\ref{thmm:romanowolf} hold with $\delta=2,\gamma=0$, and
straightforward choices of $C_{n,1},C_{n,2}, C_{n,3}$,
completing the proof.
\end{pf*}

\section*{Acknowledgements}
We are grateful to an Associate Editor and three anonymous referees for
their insightful comments.
Wei's research was supported by the National Science Foundation
(DMS-09-06568) and a career award from
NIEHS Center for Environmental Health in Northern Manhattan
(ES-009089). Zhao's research was supported
by a NIDA Grant P50-DA10075-15. The content is solely the
responsibility of the authors and does not necessarily
represent the official views of the NIDA or the NIH.

%



\printhistory

\end{document}